\documentclass{article}
\usepackage{amsmath,amssymb,graphicx,MnSymbol}
\def\MM#1{\boldsymbol{#1}}
\newcommand{\pp}[2]{\frac{\partial #1}{\partial #2}}

\def\MM#1{\boldsymbol{#1}}
\DeclareMathOperator{\diff}{d}

\DeclareMathOperator{\surf}{surf}

\usepackage{natbib, color}
\usepackage{helvet}
\usepackage{amsfonts}

\usepackage[margin=2cm]{geometry}
\newcommand{\jump}[1]{[\![#1]\!]}
\newcommand{\both}[1]{\llangle#1\rrangle}
\newcommand{\response}[1]{{\color{blue}{#1}}}

\usepackage{fancybox}
\begin{document}
\title{A compatible finite element discretisation for the
  nonhydrostatic vertical slice equations} \author{C. J. Cotter\footnote{Department of Mathematics, Imperial College London, {\texttt{colin.cotter@imperial.ac.uk}}} \hspace{1mm} and
  J. Shipton\footnote{Mathematics and Statistics, University of Exeter}} \maketitle

\begin{abstract}
  We present a compatible finite element discretisation for the
  vertical slice compressible Euler equations, at next-to-lowest order
  (i.e., the pressure space is bilinear discontinuous functions\response{)}. The
  equations are numerically integrated in time using a fully implicit
  timestepping scheme which is solved using monolithic GMRES
  preconditioned by a linesmoother. \response{The linesmoother only
    involves local operations and is thus suitable for domain
    decomposition in parallel. It allows for arbitrarily large
    timesteps but with iteration counts scaling linearly with Courant
    number in the limit of large Courant number.}  This
  \response{solver approach} is implemented using Firedrake, and the
  additive Schwarz preconditioner framework of PETSc. We demonstrate
  the robustness of the scheme using a standard set of testcases that
  may be compared with other approaches.
\end{abstract}

\section{Introduction}

This article presents numerical results for a compatible finite
element discretisation of the compressible Euler equations in a
vertical slice geometry (i.e. two dimensional with one of the
directions being the vertical). Vertical slice geometry provides an
opportunity to evaluate the performance and behaviour of numerical
discretisations for atmosphere dynamical cores using testcases that
can be run on a laptop or workstation, forming a useful step in the
development of global atmosphere dynamical cores.

The motivation for compatible finite element methods is that they
provide a naturally stable discretisation for the equations linearised
about a state of rest, i.e. no numerical stabilisation (such as
Riemann solvers, divergence damping, artificial viscosity, etc) is
required for the wave component of the solution. This stability means
that they are absent of spurious pressure modes, i.e. pressure
patterns that are highly oscillatory but lead to small numerical
gradients. Further, when the Coriolis force is introduced, they
support exactly steady geostrophic states \citep{cotter2012mixed}, and
they avoid the spurious branches of inertial oscillations
\citep{natale2016compatible} that are present in many finite element
discretisations (such as P1$_{DG}$-P2 and CR1-P0 on triangles). This
satifies the main requirements set out in
\citep{staniforth2012horizontal}, which became a wishlist for
discretisations considered in the Gung Ho project designing a new
atmosphere dynamical core for the Met Office.
\citet{natale2016compatible} also showed that if an appropriate space
is chosen for potential temperature as proposed by
\citet{guerra2016high} (guided by the Lorenz staggering for finite
difference methods), there are no spurious hydrostatic
modes. Compatible finite element methods are also flexible in allowing
to choose finite element spaces so that there are sufficient velocity
degrees of freedom per pressure degree of freedom (avoiding the
spurious inertia-gravity wave modes present in triangular C-grid
discretisations \citep{danilov2010utility}). Finally, they are
consistent on very general meshes (arbitrary triangulations with some
minimum aspect ratio, and quadrilateral meshes coming from piecewise
smooth maps applied to regular grids) and spaces can be selected of
arbitrary high order consistency. In particular, the Coriolis term is
consistent on cubed sphere meshes, avoiding an issue discovered with
the C-grid discretisation when applied to cubed sphere or dual
icosahedral meshes when formulated to satisfy the properties above
\citep{thuburn2012framework,thuburn2014mimetic}. \response{The
  numerical weather prediction community is moving towards such grids
  because they avoid the parallel bottlenecks associated with the
  poles in the latitude longitude grids in current and previous use.
  It was for this reason that compatible finite element methods were
  selected for the Gung Ho dynamical core, which is built around the
  lowest order spaces currently
  \citep{sergeev2023simulations,melvin2019mixed}.}

Compatible finite element methods also lend themselves to variational
and conservative formulations. These formulations are the result of
extension of similar schemes constructed using the C-grid finite
difference staggering
\citep{arakawa1981potential,gassmann2013global,dubos2015dynamico}.
\citet{taylor2020energy} have recently provided a formulation using
spectral element methods.

\citet{mcrae2014energy} provided a
scheme for the rotating shallow water equations that conserves energy
and enstrophy. This has been extended to bounded domains
\citep{bauer2018energy}, upwinded formulations that dissipate
enstrophy but preserve energy \citep{wimmer2020energy}, and vertical
slice model \citep{wimmer2021energy}. An alternative set of compatible
finite element spaces based on splines was presented in
\citet{eldred2019quasi}, and a related approach based on mimetic
spectral elements was presented in \citet{lee2018mixed,lee2020mixed}.

A variational discretisation (i.e., a discretisation derived from
Hamilton's principle) for the two dimensional incompressible Euler
equations was derived in \citet{natale2018variational}, which was
further developed for more general fluid models in
\citet{gawlik2020conservative,gawlik2022finite}. We are not using
variational or conservative formulations in this paper, but we note
that our formulation is rather close to them, which might be expected
to reduce spurious energy transfers between kinetic, potential and
internal energy.

The lowest order compatible finite element spaces (i.e., those spaces
with \response{RT0 for velocity and }P0 used for pressure/density) are
closely related to the C grid finite difference method popular amongst
many operational weather models \citep[for
  example]{wood2014inherently}. However, these spaces are only first
order accurate, \response{so a finite volume approach is required where one is
in effect solving for cell averaged quantities, and higher order
accuracy for those cell averaged quantities must be obtained by using
finite volume transport schemes, which have a stencil over several cells.
This was done in \citep{melvin2019mixed}.} The
advantage of using the next-to-lowest-order (NLO) spaces
(\response{RT1 is used for velocity, and }P1$_{DG}$ used for
pressure/density) is that they are naturally second order accurate, so
standard finite element formulations are sufficient. This leads to a
very clean formulation where assembly only requires to fetch data from
a single cell or from pairs of cells joined by facets. This makes it
particularly amenable to automated code generation of MPI parallel
codes using tools such as Firedrake \citep{rathgeber2016firedrake}.
\citep{shipton2018higher} presented a practical scheme for the
rotating shallow water equations using NLO spaces. In this paper we
present a practical scheme for the compressible Euler equations using
NLO spaces, evaluated using a standard suite of vertical slice model
testcases. A related scheme was coupled with a moisture model in
\citet{bendall2020compatible}, but was not benchmarked against the
tests considered here. \response{We use a fully implicit timestepping
  method, solved using a monolithic iterative solver i.e. the full
  coupled system of all variables is solved together without
  elimination. The equations are solved using GMRES, using a
  linesmoother which acts as an additive Schwarz method with each
  patch comprising the ``star'' of a vertical edge, i.e.  all degrees
  of freedom associated with the interior of the vertical column
  surrounding that edge. The columnar approach is necessary due to the
  thin domain geometry occuring in atmosphere models. Numerical
  analysis of these schemes is difficult, but they are motivated by
  finding monolithic smoothers that result in a (block) Jacobi
  iteration for the density/pressure variable after elimination within
  the patch. For example, a similar method (in local rather than
  columnar form) was proposed in
  \citet{adler2021monolithic,laakmann2022augmented}.

  There are a few reasons for
  including this approach here. First, it is useful to focus on the
  behaviour of the spatial discretisation without clouding the issue
  with aspects related to the splitting of advection and wave
  propagation, for example. In a monolithic approach, the basic code
  becomes very simple, and the challenges are exported to the
  iterative solver, which is provided here by the PETSc library.
  Second, the monolithic approach is itself new, and the paper
  provides a useful demonstration that it is useable in practice for
  numerical weather prediction.}

The rest of the article is structured as follows. Section
\ref{sec:discretisation} presents the discretisation in space and time,
together with the iterative solution strategy for the resulting implicit
system of equations, and our strategy for obtaining hydrostatic balance
reference profiles. Section \ref{sec:examples} presents the numerical
examples, and section \ref{sec:summary} provides a summary and outlook.

\section{Discretisation}
\label{sec:discretisation}
\subsection{Governing equations}
Our description of the two-dimensional \response{nonhydrostatic} dry
Euler equations can remain quite brief because our approach
discretises the equations in Cartesian coordinates (even on terrain
following coordinates or in spherical geometry).  We write the
equations in $\theta-\Pi$ (potential temperature-Exner pressure) form,
\begin{align}
  \label{eq:dudt}
  \pp{\MM{u}}{t} + 
  (\MM{u}\cdot\nabla)\MM{u} +
  f\hat{\MM{k}}\times \MM{u} + c_p\theta\nabla \Pi + g\hat{\MM{k}} & = 0, \\
  \label{eq:dthetadt}
  \pp{\theta}{t} + \MM{u}\cdot\nabla\theta & = 0, \\
  \pp{\rho}{t} + \nabla\cdot(\MM{u}\rho) & = 0, \\
  \Pi^{(1-\kappa)/\kappa} = \frac{R}{p_0}\rho\theta, & 
  \label{eq:Pi}
\end{align}
where $\MM{u}$ is the velocity, $\theta$ is the potential temperature,
$\rho$ is the density, $\Pi$ is the Exner pressure, $f$ is the
Coriolis parameter, $\hat{\MM{k}}$ is the unit vector pointing
upwards, $c_p$ is the specific heat at constant pressure, $R$ is the
gas constant, $p_0$ is the reference pressure, $g$ is the acceleration
due to gravity and $\kappa=R/c_p$. \response{Here, we have presented
  the equations as they would be read in three-dimensional form. In
  the vertical slice model, one restricts the spatial dependence of
  the fields to the $x-z$ plane, whilst retaining three Cartesian
  components of the vector field, i.e. $\MM{u}(\MM{x},t):\mathbb{R}^2\times
  [0,T]\to \mathbb{R}^3$. In a code, a simple way to implement this is
  to use a three dimensional mesh that is one cell wide in the
  $y$-direction, and periodic in that boundary condition. If a higher
  order finite element is used, it is still possible to represent
  variations in the $y$-direction. We expect $y$-independent solutions
  should be preserved by symmetry, but care should be taken to check
  that this is the case. In fact this is an advantageous approach,
  since the same code can be used for vertical slice and full 3D
  models, and validation in the vertical slice case adds to confidence
  in the full 3D model.}

In this paper, we use the vector-invariant form of the advection term,
\begin{equation}
  (\MM{u}\cdot\nabla)\MM{u} = (\nabla\times\MM{u})\times \MM{u}
  + \frac{1}{2}\left|\nabla\MM{u}\right|^2.
\end{equation}

Some of the test problems we will
consider contain a Newtonian damping term to absorb radiating internal
waves as they reach the upper boundary, requiring the addition of the following term,
\begin{equation}
\mu\hat{\MM{k}}\MM{u}\cdot\hat{\MM{k}},
\end{equation}
to the left hand side of equation \eqref{eq:dudt} where $\mu$ is the
(spatially-dependent) damping parameter specific to the test problem
(which we will specify later), and is otherwise zero.

Similarly, one of the test
problems requires the addition of viscous and diffusion terms in order
to observe convergence with mesh resolution, in which case the same
diffusion parameter $\nu$ is used. Then we add
\begin{equation}
 - \nu\nabla^2\MM{u}
\end{equation}
to the left hand side of equation \eqref{eq:dudt},
and 
\begin{equation}
 - \nu\nabla^2\theta
\end{equation}
to the left hand side of equation \eqref{eq:dthetadt}. These terms
are not necessary for stability and are only included to facilitate
convergence testing and comparison with other published results.

We call the computational domain $\Omega$, which in this paper is
always a rectangle with lateral periodic boundary conditions (or a
deformation of a rectangle to accommodate topography), and slip
boundary conditions $\MM{u}\cdot\MM{n}=0$ on the bottom and sides,
where $\MM{n}$ is the unit normal to the boundary. When diffusion and
viscosity are included, this boundary condition is augmented by
$\pp{\theta}{n}=0$ and $\pp{\MM{u}\times \MM{n}}{n}=0$.

\subsection{Spatial discretisation}
\subsubsection{Meshes and finite element spaces}
In this paper, we make use of two types of meshes. For problems with
$f=0$ and no out-of-plane component to the velocity, we construct
structured meshes of regular rectangles. For problems with $f\neq 0$,
we construct structured meshes of regular cuboids that are one cell
wide and periodic in the $y$-direction (to facilitate efficient
solution of $y$-independent problems and to allow seamless transition
to fully 3D problems). In either case, for problems with topography we
apply a terrain following transformation to the mesh of the form $z
\mapsto z + \phi(x,y,z)$, which preserves column structure (vertical
faces of cells remain vertical) but does deform rectangles into
trapezia and cuboids into trapezium prisms.

Following \citet{mcrae2016automated}, we select the finite element
spaces as follows. The velocity space $\mathbb{V}_1$ is the
Raviart-Thomas space of degree 1 (RT1) on \response{hexahedra (see
  below for more discussion about the vertical slice case)}, the
density space $\mathbb{V}_2$ is the discontinuous trilinear space
\response{(bilinear in 2D),} and temperature space
$\mathbb{V}_{\theta}$ is the tensor product of quadratic functions in
the vertical, and \response{bilinear functions (linear in 2D)} in the
horizontal. It is continuous in the vertical, and discontinuous in the
horizontal.  This choice mirrors the structure of the vertical
component of the velocity space, to facilitate the representation of
hydrostatic balance (see \citet{natale2016compatible} for an analysis
of this, and a more detailed description of these spaces). See
\citet{bercea2016structure} for information about how finite element
methods with these spaces can be used performantly.

\response{In 3D,} the $\mathbb{V}_1$ space uses a biquadratic
construction: the vertical component of velocity in the reference cube
is quadratic in the vertical direction, and linear in the horizontal
directions, with a symmetric pattern applied to the other
component(s). \response{If an $x-z$ planar mesh is used in a vertical
  slice model, then the $x-z$ components of velocity are represented
  in the RT1 space on quadrilaterals, which the $y$ component is
  discontinous and bilinear (same as the density).}  The RT1 functions
are mapped into physical cells using the Piola mapping
\citep{rognes2010efficient}. This guarantees continuity of normal
components of the functions across cell facets at the expense of
replacing polynomials by rational polynomials when the cells are
trapezia or trapezium prisms.

Since these finite element spaces all contain the complete space of
linear polynomials, classical approximation theory indicates that they
can approximate smooth functions on the reference element with second
order error in the $L^2$ norm. There is a technicality that the Piola
transformation means that $\mathbb{V}_1$ on mesh elements with terrain
following meshes do not contain all linear polynomials. However,
\citet{natale2016compatible} showed that second order approximation is
still obtained if the mapping from a rectangular domain is smooth (or,
piecewise smooth provided that the number of pieces is fixed).

The spatial discretisation of the equations
(\ref{eq:dudt}-\ref{eq:Pi}) is then to find
$(\MM{u}(t),\rho(t),\theta(t))\in
\mathring{\mathbb{V}}_1\times\mathbb{V}_2\times\mathbb{V}_{\theta}$ such that
\begin{align}
  \nonumber
  \int_\Omega \MM{w}\cdot\MM{u}_t +
   \mu\MM{w}\cdot \hat{\MM{k}}\MM{u}\cdot\hat{\MM{k}} + 
   f\MM{w}
  \cdot \hat{\MM{k}}\times \MM{u}\diff x + \MM{w}\cdot\hat{\MM{k}}g\diff x& \\
  \nonumber
  + \underbrace{
    \int_{\Omega} \nabla_h\times (\MM{w}\times \MM{u})\times \MM{u} \diff x
  - \int_{\Omega} \nabla\cdot\MM{w} \frac{1}{2}\left|\MM{u}\right|^2\diff x
  + \int_{\Gamma} \both{\MM{n}\times (\MM{u}\times\MM{w})}\cdot\tilde{\MM{u}}
  \diff S
  }_{\mbox{from }(\MM{u}\cdot\nabla)\MM{u}} & \\
  \label{eq:dudtfem}
  \underbrace{- \int_\Omega \nabla_h\cdot (\MM{w}\theta)c_p\Pi \diff x
    + \int_{\Gamma_v} \jump{\MM{n}\cdot \MM{w}\theta}c_p\{\Pi\}\diff S}_{\mbox{from }c_p\theta\nabla\Pi} & = 0, \quad \forall \MM{w}\in \mathring{\mathbb{V}}_1,\\
  \int_\Omega q\theta_t  -
  \nabla_h \cdot (\MM{u}q)\theta \diff x
  + \underbrace{\int_{\Gamma_v} \jump{q\MM{n}\cdot\MM{u}}\tilde{\theta} \diff S}_{\mbox{horizontal upwinding}}
  + \underbrace{\int_{\Gamma} C_0 h^2 |\MM{u}\cdot\MM{n}|
  \jump{\nabla_h q}\cdot\jump{\nabla_h\theta}\diff S}_{\mbox{edge stabilisation}} & = 0,
  \quad \forall q \in \mathbb{V}_{\theta}, \\
  \int_\Omega \phi \rho_t  - \nabla_h\phi \cdot \MM{u}\rho \diff x
  + \underbrace{\int_\Gamma \jump{\phi\MM{u}\cdot\MM{n}}\tilde{\rho}\diff S}_{\mbox{upwinding}} &
  = 0, \quad \forall \phi\in\mathbb{V}_2,
\end{align}
where 
\begin{enumerate}
  \item $\nabla_h$ indicates the ``broken'' gradient obtained by
    evaluating the gradient separately in each cell,
  \item $\Gamma$ is the set of interior facets (with $\Gamma_v$ the
    set of vertical facets between adjacent columns),
  \item $\jump{\psi}$ is the ``jump'' operator applied to a quantity
    $\psi$ returning $\psi^+-\psi^-$ where each interior facet has
    sides arbitrarily labelled ``+'' and ``-'' (and noting that
    $\MM{n}$ is the unit normal to a facet with $\MM{n}^+$ pointing
    into the $-$ side and \emph{vice versa} for $\MM{n}^-$),
  \item $\tilde{\psi}$ indicates the upwind value of any quantity
    $\psi$ on the facet i.e. the value on the side where
    $\MM{u}\cdot\MM{n}\geq 0$ (making an arbitrary choice not
    affecting the result when $\MM{u}\cdot\MM{n}=0$),
  \item $\{\psi\}$ is the average operator returning
    $(\psi^++\psi^-)/2$ when evaluated on a facet,
  \item $\Pi$ is defined according to \eqref{eq:Pi} but now applied to
    the numerical approximations $\rho$ and $\theta$ (so that we do
    not separately solve for an independent variable $\Pi$),
  \item $h$ is \response{an} estimate of the cross facet meshscale defined by
    $\{V_e\}/A_f$ with $V_e$ being the (cellwise constant) cell
    volume, and $A_f$ being the (facetwise constant) facet area,
  \item $C_0$ is an edge stabilisation constant (chosen here to be
    $2^{-7/2}$ following the suggested scaling $p^{-7/2}$ with
    polynomial degree $p$ given by \citet{burman2007continuous}), and
  \item $\mathring{\mathbb{V}}_1$ is the subspace of $\mathbb{V}_1$
    containing all functions that satisfy the $\MM{u}\cdot\MM{n}=0$
    boundary condition on the top and bottom of the domain.
\end{enumerate}
The treatment of the velocity advection term was first demonstrated
(for a incompressible Boussinesq vertical slice model) by
\citet{yamazaki2017vertical}, inspired by
\citet{natale2018variational} but different from the energy conserving
form for compressible Euler equations proposed in
\citet{wimmer2021energy}. The surface term deals with the fact that
the velocity space does not have continuous tangential components in
general; it also has a stabilising effect as examined in
\citet{natale2017scale}. The treatment of the pressure gradient term
first appeared in \citet{natale2016compatible}, and has been used in a
modified energy conserving form in \citet{wimmer2021energy}, and with
lowest order spaces 
in \citet{bendall2020compatible} and as part of the Gung Ho dynamical
core formulation in \citet{melvin2019mixed}. The treatment of the
potential temperature advection term deviates here from the ``SUPG''
formulation proposed in \citet{yamazaki2017vertical}, using an edge
stabilisation proposed by \citet{burman2005unified} instead, in
combination with standard discontinuous Galerkin style upwinding on
vertical faces. Even though the temperature space is only continuous
in the vertical direction, we found that edge stabilisation was
necessary on all faces to achieve stable results.

It should be noted that in the case of spatially varying topography,
it is not possible to compute integrals exactly because of the
$\det(J)$ in the reciprocal appearing in the velocity basis functions
(due to the Piola mapping), where $J$ is the Jacobian of the reference
element to mesh element mapping. Thus, we have to approximate the
integrals using numerical quadrature, and we have to select a
quadrature degree for each term so that the discretisation is stable
and sufficiently accurate. Here we just take the suggested quadrature
degree from the heuristic computed in UFL \citep{alnaes2014unified}
which selects a (rather generous) quadrature rule for approximate
integration, which is certainly sufficient for stability and second
order consistency.  In fact, many of the integrals are nevertheless
computed exactly, due to fortuitious cancellations that are discussed
in \citep{cotter2014finite}.  Another source of inexact quadrature is
the function $\Pi$ which involves a noninteger power of $\rho$ and
$\theta$.

When test cases require a viscous term, we use the symmetric interior
penalty discretisation (to deal with discontinuities in the tangential
component of the velocity) as described in \citep{cockburn2007note},
in which case the term
\begin{equation}
  \nu\left(\int_\Omega \nabla_h \MM{u}:\nabla_h \MM{v} \diff x -
  \int_\Gamma \jump{\MM{v}\otimes\MM{n}}:\{\nabla \MM{u}\}\diff S -
  \int_\Gamma \jump{\MM{u}\otimes\MM{n}}:\{\nabla \MM{v}\}\diff S +
  \int_\Gamma \frac{\eta}{h}
  \jump{\MM{u}\otimes\MM{n}}:\jump{\MM{v}\otimes\MM{n}}\diff S\right)
\end{equation}
is added to the left hand side of Equation \eqref{eq:dudtfem}, where
$\nu$ is the dynamic viscosity, $:$ indicates the double contraction
for tensors (so that $A:B=\sum_{ij}A_{ij}B_{ij}$), and $\eta$ is a
dimensionless penalty parameter, chosen here to have the value 10
which is experimentally determined to produce a stable discretisation
for both $\mathbb{V}_1$ and $\mathbb{V}_{\theta}$. A similar formula
(with the same parameters, but adapted to scalar fields) is used when
test cases require a diffusion term in the potential temperature
equation.

\subsection{Time discretisation and iterative solver}
The time discretisation used is the implicit midpoint rule, a fully
implicit second order method. This is obtained by replacing time
derivatives by time differences e.g. $\theta_t \mapsto (\theta^{n+1}-
\theta^n)/\Delta t$, and replacing all other occurrences of fields by
their midpoint values e.g. $\theta \mapsto
(\theta^{n+1}+\theta^n)/2$. This nonlinear system is then solved for
the $n+1$ variables using the linesearch Newton method provided by
PETSc \citep{balay2020petsc} using Firedrake's
NonlinearVariationalSolver object. The resulting Jacobian systems are
solved using GMRES applied to the full monolithic
velocity-density-temperature system, preconditioned by an additive
Schwarz method, which does direct solves over column patches
surrounding one vertex on the base mesh (in vertical slice models this
corresponds to forming a patch out of two neighbouring cells). In
fact, a ``star'' patch is used, which neglects degrees of freedom
associated with the horizontal boundaries of the patch. \response{This
  patch is formed out of the ``star iteration'' (see Section 4 in
  \citep{farrell2019augmented}) applied to a vertex in the base mesh,
  before extruding up the column.} These direct solves couple all
three fields, and the direct solve uses PETSc's own LU factorisation
algorithm using a reverse Cuthill-McGee ordering to reduce the fill
in. This algorithm can be and is parallelised using domain
decomposition, which was automated in our implementation using
Firedrake. In all of our test cases, Newton's method converges in 2-3
iterations, and GMRES converges in 10-35 iterations depending on the
flow complexity, independent of mesh size provided that a constant
Courant number is maintained whilst decreasing the mesh size.

\response{This solver approach is necessary because we chose to use a
  fully implicit timestepping scheme, solved using Newton's method.
  This requires us to be able to solve systems $Jx=b$ where $J$ is the
  full Jacobian obtained by linearisation about the state at each
  Newton iteration. In the case where $J$ is the linearisation about a
  state of rest, $J$ can be reduced to a positive definite Schur
  complement using the hybridisation technique provided the
  discretisation proposed in \citet{betteridge2022hybridised} is used.
  However, this technique does not work when the linearisation is
  calculated for states with nonzero velocities, since these
  introduce additional intercell coupling.}

\subsection{Hydrostatic balance}
The test cases that we consider here require the computation of a
background density profile $\rho_b$ that is in hydrostatic balance
given the specified potential temperature $\theta_b$. \response{We use
  a nonhydrostatic model, but testcases require a initialisation at a
  state of hydrostatic balance, which we describe here.}  To avoid
unphysical motion at lower resolutions we compute this hydrostatic
balance numerically, i.e. we require that
\begin{equation}
  \label{eq:hydros}
  - \int_\Omega \nabla_h\cdot (\MM{w}\theta_b)c_p\Pi_b \diff x
  + \int_\Gamma \jump{\MM{n}\cdot \MM{w}\theta_b}c_p\{\Pi_b\}\diff S
  + \int_\Omega g\MM{w}\cdot\hat{\MM{k}}\diff x = 0, \quad
  \forall \MM{w}\in \mathring{\mathbb{V}}_{1,v},
\end{equation}
where $\mathbb{V}_{1,v}$ the vertical subspace of $\mathbb{V}_1$,
and $\mathring{\mathbb{V}}_{1,v}$ is the corresponding vertical subspace of
$\mathring{\mathbb{V}}_1$.\footnote{\response{In the presence of orography, $\mathbb{V}_{1,v}$ does not contain all of the vertical component of the velocity, and solving
  this equation will lead to pressure gradient errors. This is described further
  in \citet{natale2016compatible}.}}
  We note that $\Pi_b$ is just a local
  nonlinear function of $\rho_b$ and $\theta_b$. As discussed in
  \citet{natale2016compatible}, we can solve this equation for $\rho_b$
  by introducing an auxiliary variable $\MM{v}$ (which will turn out
  to vanish), and solving the coupled system
\begin{align}
  \label{eq:square hydros 1}
  \int_{\Omega}\MM{w}\cdot\MM{v}\diff x - \int_\Omega \nabla_h\cdot
  (\MM{w}\theta_b)c_p\Pi_b \diff x + \int_\Omega
  g\MM{w}\cdot\hat{\MM{k}}\diff x + \int_{\partial \Omega_0} c_p\MM{w}\cdot\MM{n}
  \theta_b\Pi_0\diff S
  & = 0, \quad \forall \MM{w}\in \tilde{\mathbb{V}}_{1,v},
  \\
  \label{eq:square hydros 2}
  \int_\Omega \nabla_h\cdot
  (\MM{v}\theta_b)c_p\phi \diff x & = 0, \quad \forall \phi\in \mathbb{V}_2,
\end{align}
where $\partial \Omega_0$ is the surface where we have chosen to set
$\Pi_b=\Pi_0$ as a boundary condition (the bottom boundary for the
test cases in this paper), $\tilde{\mathbb{V}}_{1,v}$ is the subspace
of $\mathbb{V}_{1,v}$ containing functions that satisfy
$\MM{w}\cdot\MM{n}=0$ on the opposite boundary (the top boundary for
the test cases in this paper). We note that the surface integral in
Equation \eqref{eq:hydros} vanishes when $\MM{w}\in \mathbb{V}_{1,v}$
because $\MM{w}\cdot\MM{n}=0$ on vertical faces, and we note that the
weak boundary condition integral in \eqref{eq:square hydros 1} over
$\partial\Omega_0$ vanishes when $\MM{w}\in
\mathring{\mathbb{V}}_{1,v}\subset\tilde{\mathbb{V}}_{1,v}$. Since
$\MM{v}$ vanishes at the solution, we recover the hydrostatic
condition \eqref{eq:hydros}.  The system (\ref{eq:square hydros
  1}-\ref{eq:square hydros 2}) decouples into independent columns,
which we solve using Newton's method. \citet{natale2016compatible}
showed that the linearisation around a given state of this system is
well-posed, and we solve the resulting linear systems for
$(\rho_b,\MM{v})$ updates directly. We find an initial guess for
$\rho_b$ by first solving the corresponding linear system for
$(\Pi_b,\MM{v})$ where $\Pi_b\in \mathbb{V}_2$ is taken as an
independent variable instead of a local function of $\rho_b$ and
$\theta_b$.  We then project the formula for $\rho_b$ in terms of
$\theta_b$ and $\Pi_b$ into $\mathbb{V}_2$. If we use the result as an
initial guess then Newton's method converges in 2-3
iterations. \response{Note that this computation is only done to
  compute a set of initial conditions, and is not used in the forward
  model which is always nonhydrostatic.}

\section{Computational examples}
\label{sec:examples}
In this section we demonstrate our discretisation approach using the
suite of test problems considered in
\citet{melvin2010inherently}. This suite tests the vertical slice
discretiation on a range of flows that are relevant to numerical
weather prediction, including acoustic and gravity waves and flows
driven by both buoyancy and orography. The gravity wave and orographic
wave tests are run in both the hydrostatic and nonhydrostatic regimes.
To avoid the need to refer back, we provide a brief summary of the
test problems here. Some constants that are consistent across all
tests are provided in Table \ref{tab:constants}.
\begin{table}
  \begin{center}
  \begin{tabular}{|c|c|c|c|c|c|c|}
    \hline
    $g$ & $N$ & $f$ & $c_p$ & $R$ & $p_0$ & $c_v$ \\
    \hline
    9.810616$ms^{-2}$ & $10^{-2}s^{-1}$ & $10^{-4}s^{-1}$ &
    1004.5$J kg^{-1}K^{-1}$ &
    287 $J kg^{-1}K^{-1}$ & 1000 $hPa$ &
    717 $J kg^{-1}K^{-1}$ \\
    \hline
  \end{tabular}
  \end{center}
  \caption{\label{tab:constants} Some constants that take consistent
    values across all the computational examples. The exceptions are
    the nonhydrostatic variants of the gravity wave and orographic
    flows and the density current which are non-rotating so $f=0$; the
    density current also has an isentropic background state so $N=0$.}
  \end{table}
In each test we have used the same timestep $\Delta t$ values as
\citet{melvin2010inherently}, and double the values of $\Delta x$ and
$\Delta z$, which ensures the same number of degrees of freedom
(because we are using NLO spaces). A summary of \response{resolutions,
  timesteps, number of iterations, and the timings} is provided in
Table \ref{tab:stats}.
\begin{table}
  \begin{center}
  \begin{tabular}{|c|c|c|c|c|c|c|}
    \hline
    Testcase & Columns & layers & $\Delta t$ & GMRES its per step & Time & Cores \\
    \hline
    Nonhydrostatic gravity wave & 150 & 5 & 12s & 31.512   & 48.789s & 16 \\ \hline
    Hydrostatic gravity wave & 300 & 10 & 100s & 31.022   & 227.51s & 16 \\ \hline
    Nonhydrostatic mountain & 180 & 70 & 5s &  27.239 & 1912.5 & 16 \\ \hline
    Hydrostatic mountain & 100 & 60 & 20s & 33.959   & 6521.6s & 16 \\ \hline
    density current $\Delta x=800m$ & 64 & 8 & 4s & 62.8  & 60.04s & 16 \\  \hline
    density current $\Delta x=400m$ & 128  & 16 & 2s & 55.05    & 208.4s & 16 \\ \hline
    density current $\Delta x=200m$ & 256 & 32 & 1s & 49.11   & 1157s & 16 \\ \hline
    density current $\Delta x=100m$ & 512 & 64 & 0.5s & 43.46   & 8452s & 32 \\ \hline
    Sch\"ar mountain &  100 & 50 & 8s & 49.802   & 1699.1s & 16 \\  \hline
    Sch\"ar mountain & 100  & 50 & 40s & 274.756   & 1342.6s & 16 \\ \hline
  \end{tabular}
  \end{center}
  \caption{\label{tab:stats} Table showing \response{the resolution, $\Delta t$},
    number of GMRES iterations per timestep (over all the Newton
    iterations), the time taken and the number of cores used, for each
    of the testcases. The timings include the nonhydrostatic
    initialisation for the density and solver setups.}
\end{table}

\subsection{Gravity waves}
\label{sec:gw tests}
This test case is the ``Hello World!'' of vertical slice test
problems, first proposed in \citet{skamarock1994efficiency}. There are
two versions of the test case, the nonhydrostatic flow regime version
with velocity constrained to the $x-z$ plane and consequently $f=0$,
and the hydrostatic flow regime version with 3D velocity and
$f=10^{-4}s^{-1}$. \response{Note that we are solving the nonhydrostatic
  equations in both versions, and the naming just describes the asymptotic
  flow regime and not the equations solved.}
The domain is given by $L/2\leq x \leq L/2$ and $0
\leq z \leq H$ where $L=3\times 10^5m$ in the nonhydrostatic case and
$L=6\times 10^6m$ in the hydrostatic case. In both cases the height is
$H=10^4m$ and there are periodic boundary conditions in the horizontal
direction.

In both cases, the potential temperature is initialised to a
background profile
\begin{equation}
  \theta_b = T_{\surf}\exp(N^2z/g),\quad T_{\surf}=300K,
\end{equation}
before solving for the hydrostatically balanced $\rho_b$. Then a
perturbation is added to the potential temperature, 
\begin{equation}
  \Delta \theta = \Delta \theta_0\frac{\sin(\pi z/H)}{1 + x^2/a^2},
\end{equation}
where $\Delta \theta_0=10^{-2}K$ and $a=5\times 10^3m$ for the
nonhydrostatic flow regime and $a=10^5m$ for the hydrostatic flow
regime. In both cases, the horizontal velocity in the $x$-direction
is
initialised to $20ms^{-1}$ and the other components to zero. In the
hydrostatic case, an additional forcing term is introduced to balance
the Coriolis force, adding $f\times(0,-20,0)$ to the left hand side of
Equation \eqref{eq:dudt}.

\begin{figure}
    \centerline{\includegraphics[width=18cm]{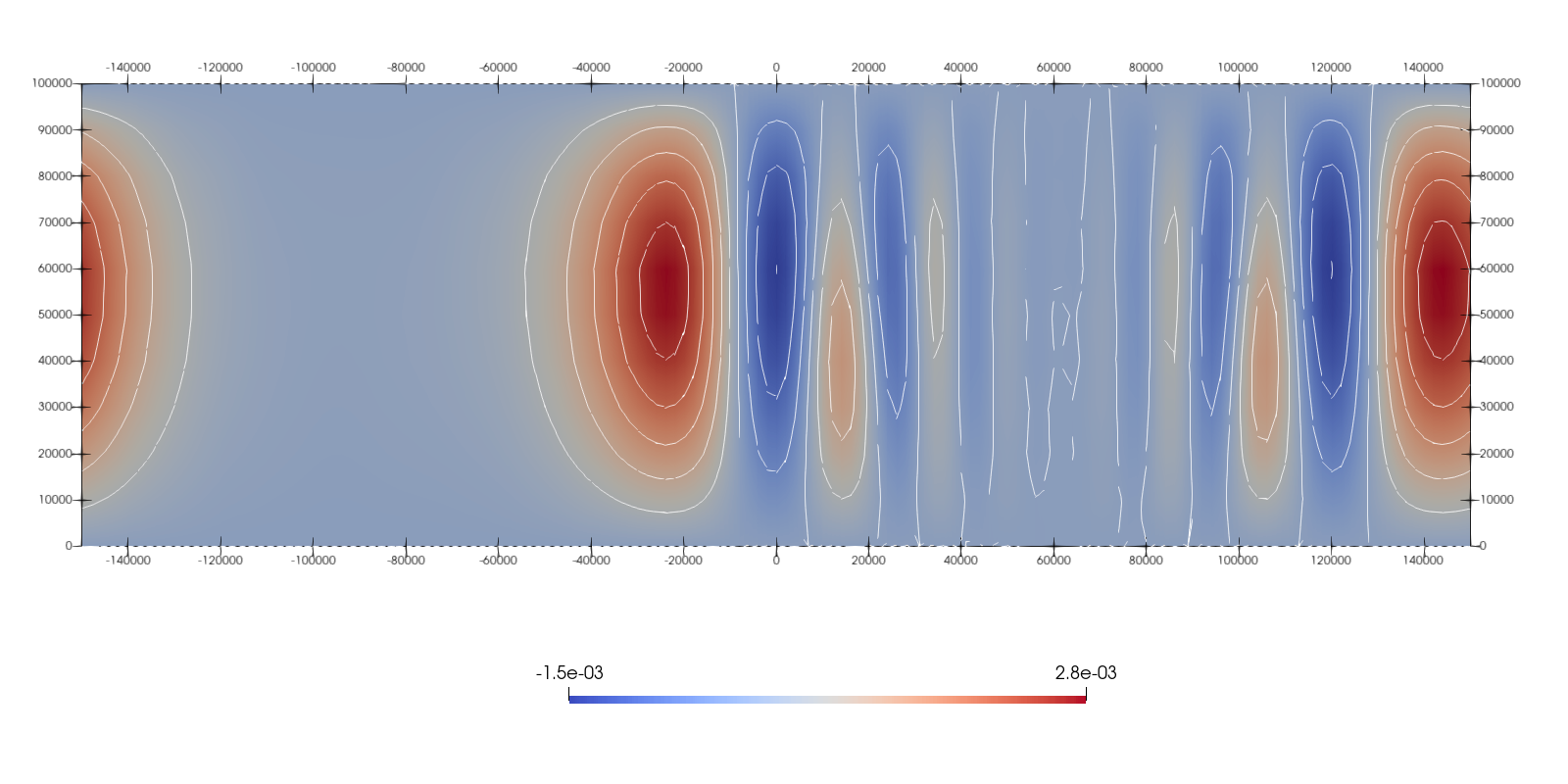}}
    \caption{\label{fig:gw nh}Contour plot of vertical velocity for the
      nonhydrostatic gravity wave at time $t=3000s$ at resolution
      $\Delta x=2000m$, $\Delta z=2000m$, and $\Delta t=12s$. Contours
      are drawn every $5\times 10^{-4}ms^{-1}$.}
\end{figure}

\begin{figure}
  \centerline{\includegraphics[width=22cm]{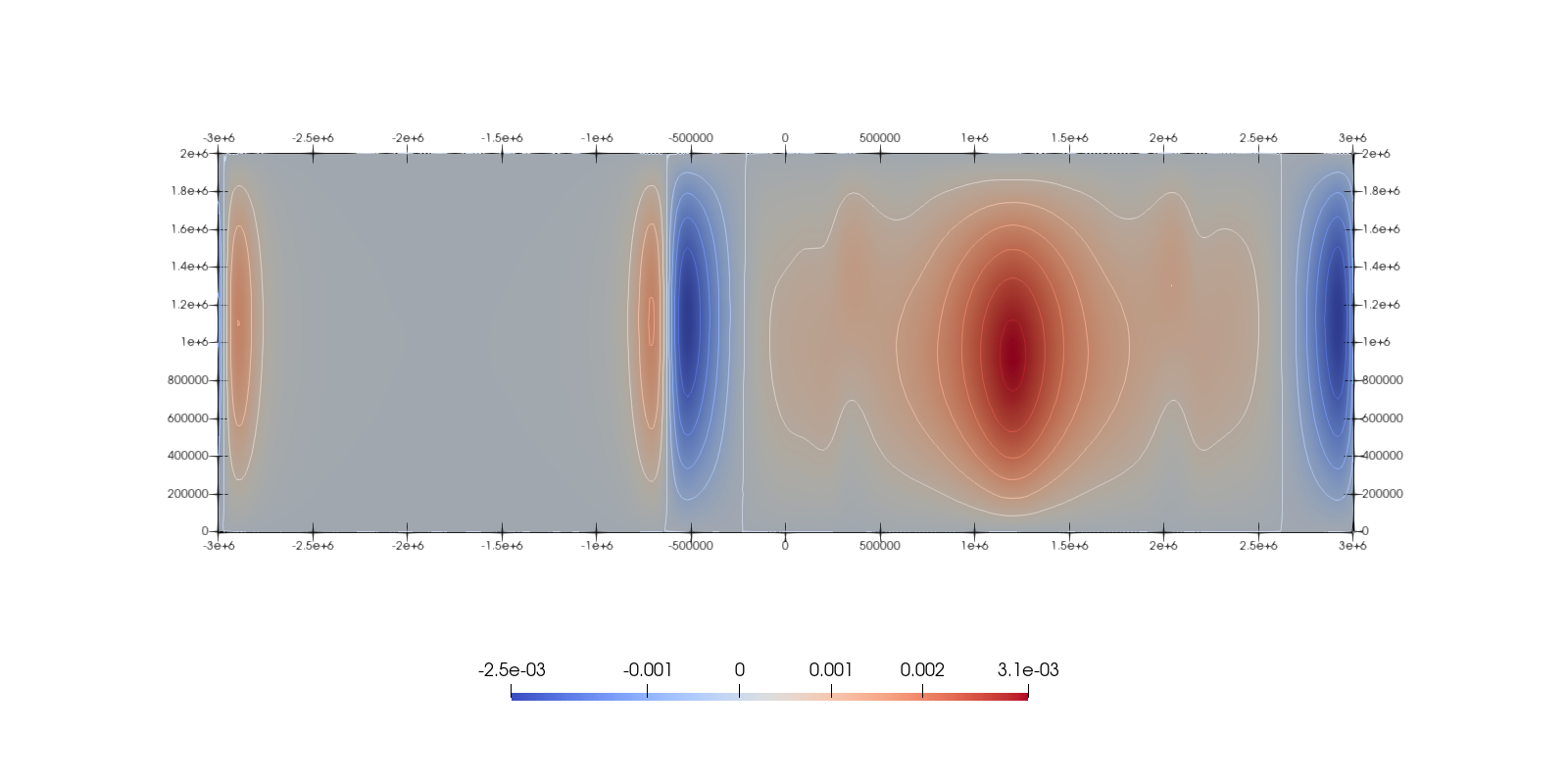}}
  \caption{\label{fig:gw h}Contour plot of vertical velocity for the
    hydrostatic gravity wave at time $t=60000s$ at resolution $\Delta
    x=20000m$, $\Delta z=1000m$, and $\Delta t=100s$. Contours are
    drawn every $5\times 10^{-4}ms^{-1}$.}
\end{figure}

Plots of the nonhydrostatic and hydrostatic flow regime solutions are
shown in Figures \ref{fig:gw nh} and \ref{fig:gw h}, respectively.
These solutions closely match the results from
\citet{melvin2010inherently} at similar resolutions.

\subsection{Density current}
This test is taken from the classic intercomparison project of
\citet{straka1993numerical}, simulating a dense bubble in an
isentropic, hydrostatic atmosphere. The domain is
$-L/2\leq x \leq L/2$ where $L=51200m$, and $0 \leq z \leq  H=6400m$ with
periodic boundary conditions in the horizontal direction.

The background temperature profile is chosen to be isentropic,
i.e. the potential temperature is constant. In this case the
background potential temperature is $\theta_b=300K$. The background
density profile is then obtained by solving for hydrostatic balance
for this potential temperature profile, with the boundary condition
$\Pi=1$ on the bottom boundary. We then apply a perturbation
to the temperature, 
\begin{equation}
  \Delta T = \left\{
  \begin{array}{cr}
    0 & \mbox{if }L_r > 1, \\
    -15\left(\cos(\pi L_r) + 1\right)/2 & \mbox{otherwise},
  \end{array}\right.
\end{equation}
at constant pressure $p$, where
\begin{equation}
  L_r = \sqrt{\left(\frac{x \response{- x_c}}{x_r}\right)^2
+    \left(\frac{z-z_c}{z_r}\right)^2},
\end{equation}
and $(x_c,x_r)=(0m, 4000m)$, $(z_c,z_r)=(3000m,2000m)$.  This
corresponds to making the potential temperature perturbation $\Delta
\theta =\Delta T/\Pi$, where $\Pi$ is the background Exner pressure
profile, and then perturbing the density according to $(\rho+\Delta
\rho)(\theta+\Delta \theta) = \rho\theta$. In our implementation,
$\Delta \theta$ was computed by the following steps:
\begin{enumerate}
\item Project the formula $\Delta T/\Pi(\rho_b,\theta_b)$ into
  $\mathbb{V}_{\theta}$, where $\rho_b\in\mathbb{V}_2$ and
  $\theta_b\in\mathbb{V}_{\theta}$ are the previously computed initial
  conditions for density and temperature respectively.  This is the
  initial condition $\theta_0$ for potential temperature.
\item Project the formula $\rho_b\theta_0/\theta_b$ in to $\mathbb{V}_2$.
  This is the initial condition $\rho_0$ for density.
\end{enumerate}
These projections are $L^2$ projections using incomplete but high order
quadrature for the integrals over the formulae, similar to those used
in the dynamical equations.

Without viscosity, this problem becomes ill-posed after forming a
singular vorticity structure in finite time.  Hence, \response{to be
  able to compare between models}, a viscosity of $\nu=75m^2s^{-1}$ is
included in the velocity equation and a diffusivity of
$\kappa=75m^2s^{-1}$ is included in the potential temperature
equation.

For this testcase, the comparison is made after 15 minutes. Contour
plots at various resolutions are shown in Figure \ref{fig:straka
  contours}, and some summary statistics are provided in Table
\ref{tab:straka}. There is quite a bit of variation between models for
the precise solution at this point, as there is a lot of small scale
structure forming, but our contour plots look broadly similar to those
of \citet{melvin2010inherently}. In particular, the density current has
reached a very similar location, estimated as the maximum $x$
coordinate over all cells where $\Delta\theta<0$. We see that the
dissipation of the minima of $\Delta \theta$ is much weaker at the
higher resolutions, and although there is still a substantial
overshoot past the initial maxima of 0, this may be because the
solution produces finer filaments of potential temperature at higher
resolutions, which is more challenging for the advection scheme.

\begin{figure}
  \centerline{\includegraphics[width=8cm]{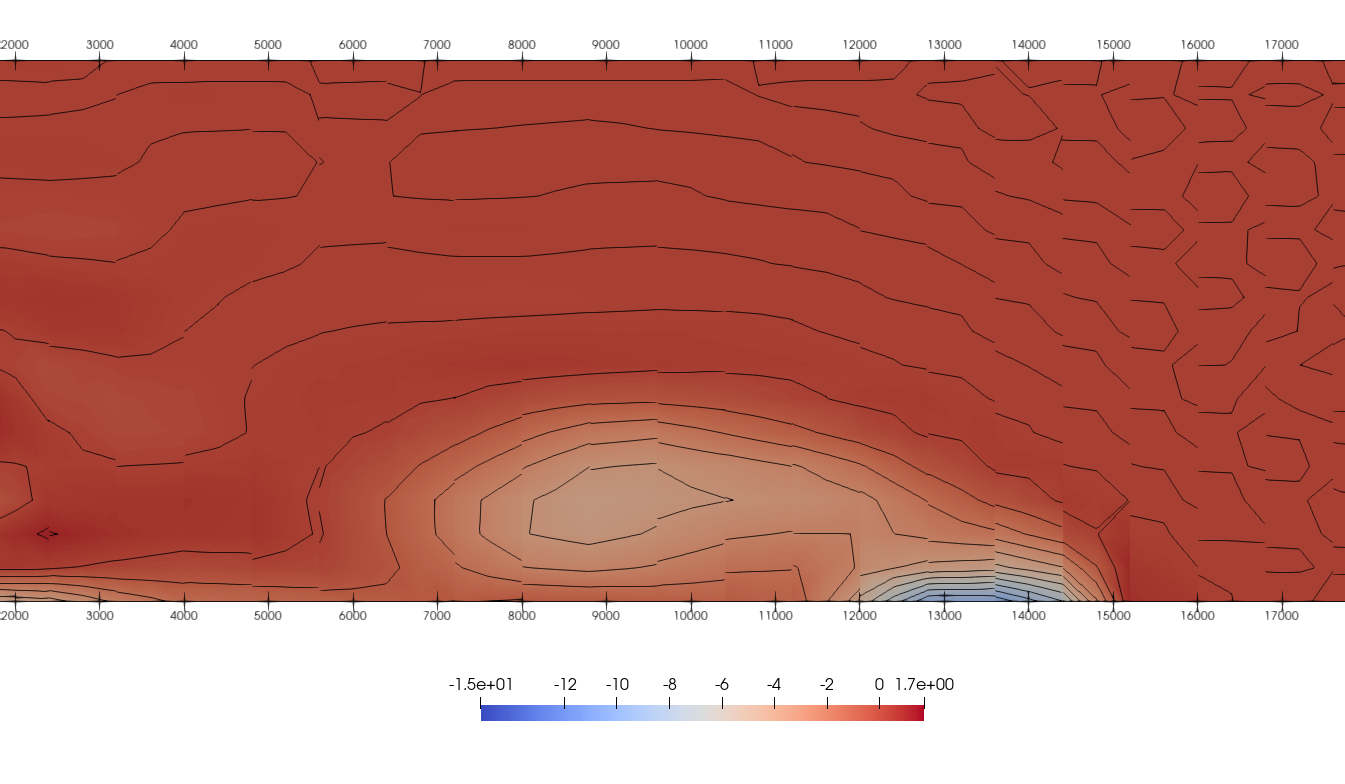}
    \includegraphics[width=8cm]{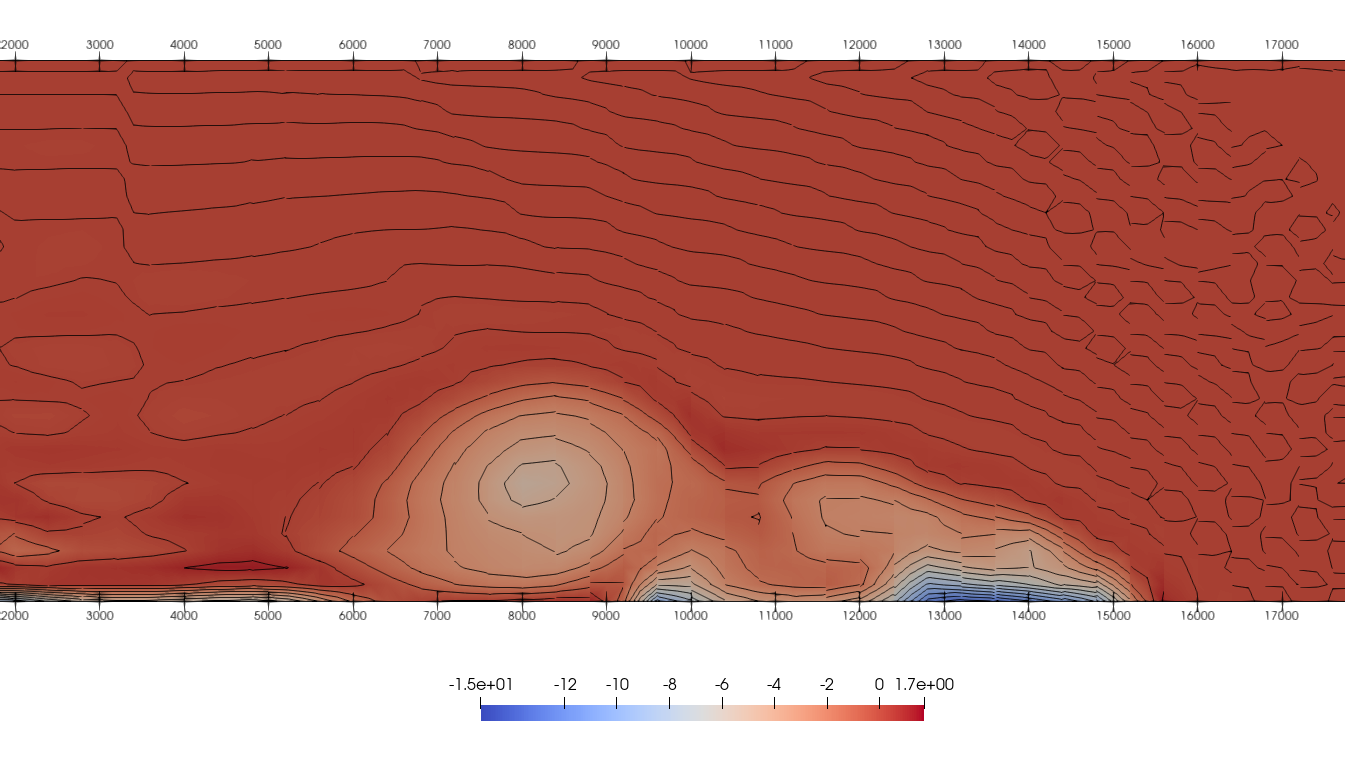}}
  \centerline{\includegraphics[width=8cm]{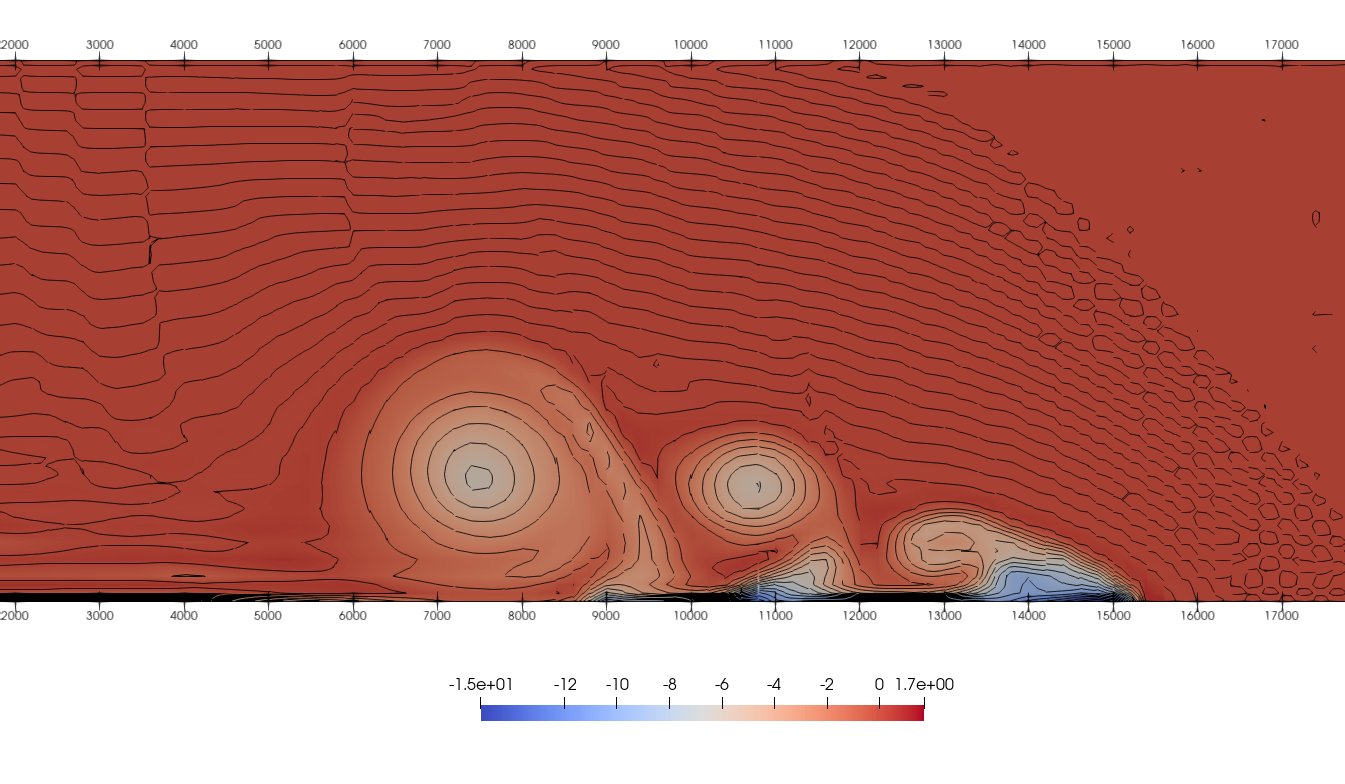}
    \includegraphics[width=8cm]{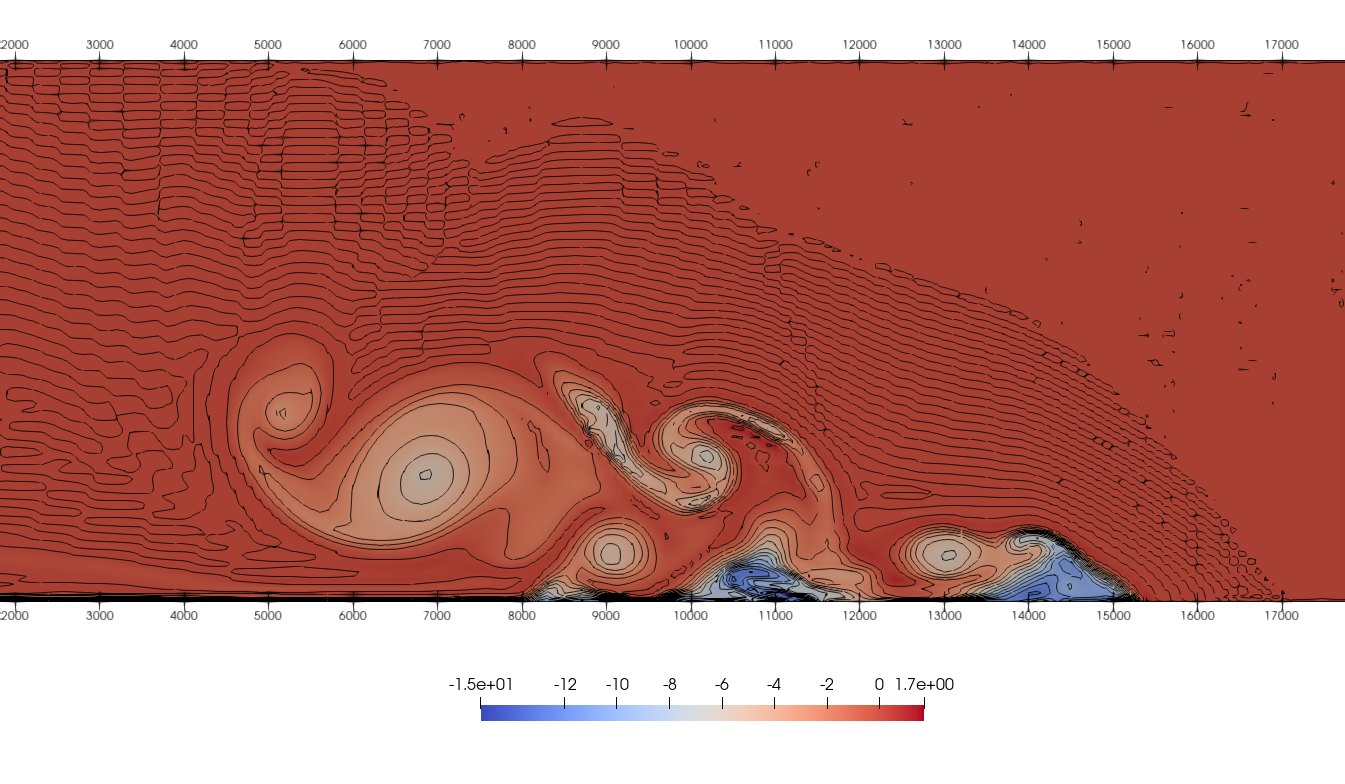}}
  \caption{\label{fig:straka contours}For the density current
    testcase: contours of $\Delta \theta$ at 15 minutes. Contour
    intervals are 1$K$, with the first contour being at -1 $K$. From
    left to right and top to bottom, the resolutions are $(\Delta x,
    \Delta t)=(800m,4s)$, $(400m, 2s)$, $(200m, 1s)$ and $(100m,
    0.5s)$ respectively.}
\end{figure}

\begin{table}
  \begin{center}
\begin{tabular}{|c|c|c|c|c|}
  \hline
  $\Delta x$ & $\Delta t$ & $\Delta\theta_{max}$ & $\Delta\theta_{min}$ & front location \\
  \hline
  800m & 4s & 1.3448 & -11.6445 & 14800 \\ 
  400m & 2s & 1.2293 & -13.4614 & 15000 \\
  200m & 1s & 0.9023 & -16.3183 & 15300 \\
  100m & 0.5s & 1.4226 & -15.0101 & 15250 \\
  \hline
\end{tabular}
\end{center}
\caption{\label{tab:straka}For the density current testcase: minimum,
  and maximum values of $\Delta\theta=\theta-\theta_b$ measured at 15
  minutes, together with the front location (estimated as the maximum $x$
  coordinate over all cells where $\Delta\theta<0$). The image has been cropped to focus on the right-propagating current.}
\end{table}

\subsection{Flow over a mountain}

This test problem simulating small amplitude lee waves generated by
flow over a (small) mountain also has two versions, the nonhydrostatic
flow regime version with velocity constrained to the $x-z$ plane and
consequently $f=0$, and the hydrostatic flow regime version with 3D
velocity and $f=10^{-4}s^{-1}$. \response{Note that we are solving the
  nonhydrostatic equations in both versions, and the naming just
  describes the asymptotic flow regime and not the equations solved.}
The domain is given by $L/2\leq x \leq L/2$ and $0 \leq z \leq H$
where $L=144000m$ and $H=35000m$ in the nonhydrostatic case and
$L=240000m$ and $H=50000m$ in the hydrostatic case. In both cases
there are periodic boundary conditions in the horizontal direction.

The mountain has a ``Mount Agnesi'' profile, with the bottom boundary
moved to $z_s(x) = \frac{a^2}{x^2+a^2}$, giving a mountain of height
$1m$. In our implementation, we use a simple terrain following mesh
with the rectangular domain transformed according to $(x,y,z)\mapsto
(x,y,z + z_s(H-z)/H)$. In the nonhydrostatic flow regime, $a=10000m$,
and in the hydrostatic regime, $a=1000m$.

In the nonhydrostatic flow regime, a stratified background flow is
initialised according to the description of Section \ref{sec:gw
  tests}, with $T_{surf}=300K$. In the hydrostatic flow regime, the
stratification is isothermal, i.e. constant temperature at
$T=T_{surf}=250K$. This does still imply a varying potential
temperature, with profile
\begin{equation}
  \theta = T_{surf}\exp\left(\frac{gz}{T_{surf}c_p}\right).
\end{equation}
In both cases, the density is initialised by solving numerically for a
hydrostatic profile, with boundary condition $\Pi=1$ at $z=0$.  This
requires additional calculation because with the topography, the
bottom boundary is not at $z=0$ everywhere. To address this, we
calculate the boundary condition for $\Pi$ at the top of the domain
that produces the value $\Pi=1$ at $z=0$, using a value on the bottom
of the domain away from the mountain.

The nonhydrostatic test is initialised with a horizontal velocity
$u=10ms^{-1}$ and the hydrostatic test is initialised with a
horizontal velocity $u=20ms^{-1}$.  As for the gravity wave test, in
the hydrostatic case, an additional forcing term is introduced to
balance the Coriolis force, adding $f\times(0,-20,0)$ to the left hand
side of Equation \eqref{eq:dudt}. In both cases, the initial velocity is not
compatible with the boundary condition at the mountain, so this causes
a pressure wave propagating at ground level, radiating waves that
interfere with the stationary lee wave pattern that accumulates over
time. This is often referred to as a test of robustness of the
discretisation.

To prevent the lee waves \response{(and the initial pressure waves
  caused by the sudden appearance of the mountain at time 0)}
reflecting off the top boundary, an absorbing term in the vertical
velocity is added in the top layer, with profile
\begin{equation}
  \label{eq:absorb}
  \mu(z) = \left\{ \begin{array}{rl}
    0, & z < z_B, \\
    \bar{\mu} \sin^2\left(\frac{\pi}{2}\left(\frac{z-z_B}{H-z_B}\right)\right),
    & z \geq z_B,
  \end{array}
  \right.
\end{equation}
where $\bar{\mu}$ is a constant and $z_B$ is the height of the bottom
of the absorbing layer. For the hydrostatic test, $z_B=H-2\times
10^4m=3\times 10^4m$ and $\bar{\mu}\Delta t=0.3s$. For the nonhydrostatic
test, $z_B=H-10^4m=2.5\times 10^4m$ and $\bar{\mu}\Delta t=0.15s$.

\begin{figure}
  \centerline{\includegraphics[width=18cm]{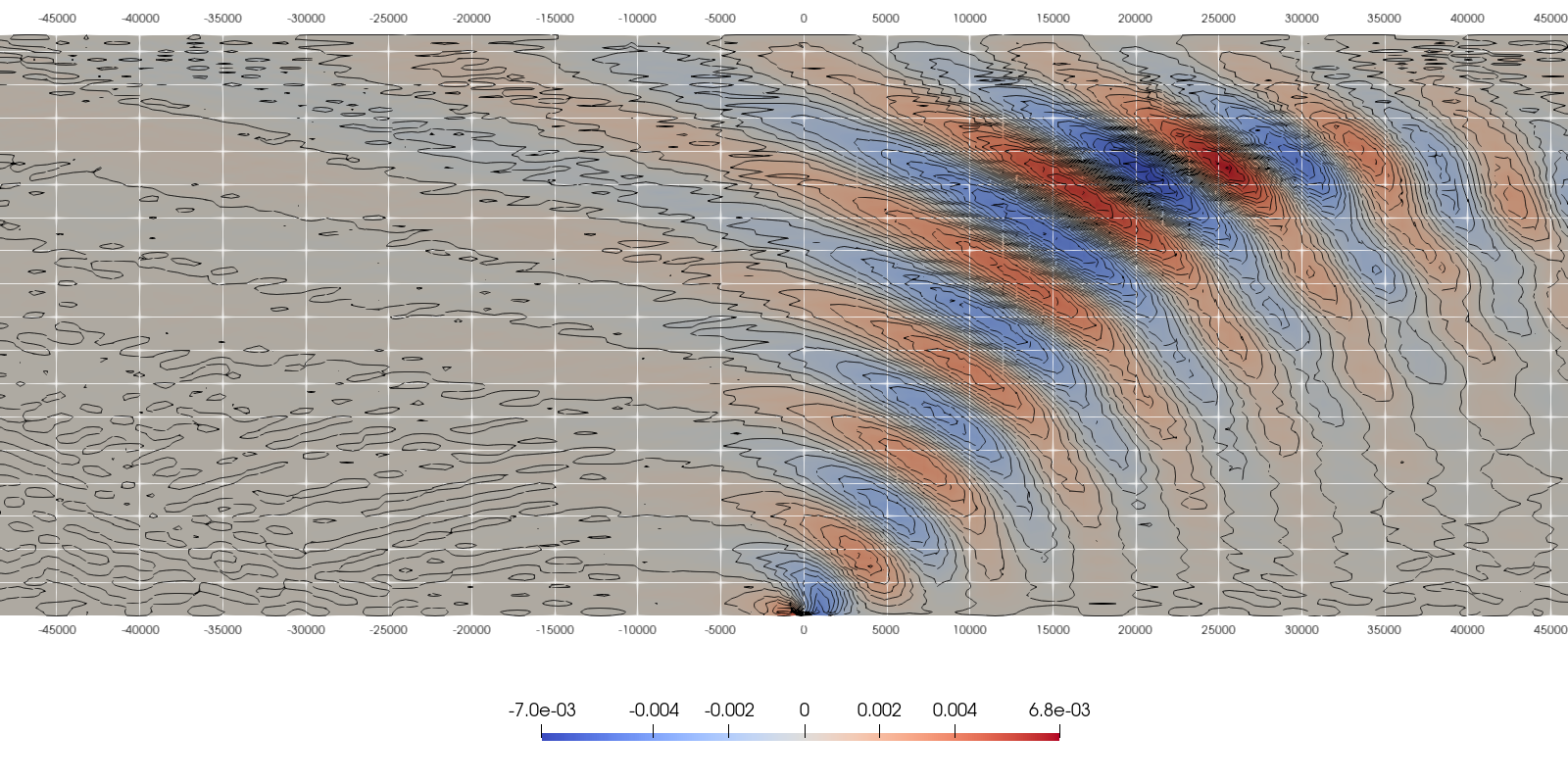}}
    \caption{\label{fig:ag nh}Contours of
  vertical velocity for the nonhydrostatic mountain wave test at
  $t=9\times 10^3s$. Contour intervals are every $5\times
  10^{-4}ms^{-1}$\response{; horizontal lines are plotted every 2000m.}}
\end{figure}

\begin{figure}
  \centerline{\includegraphics[width=18cm]{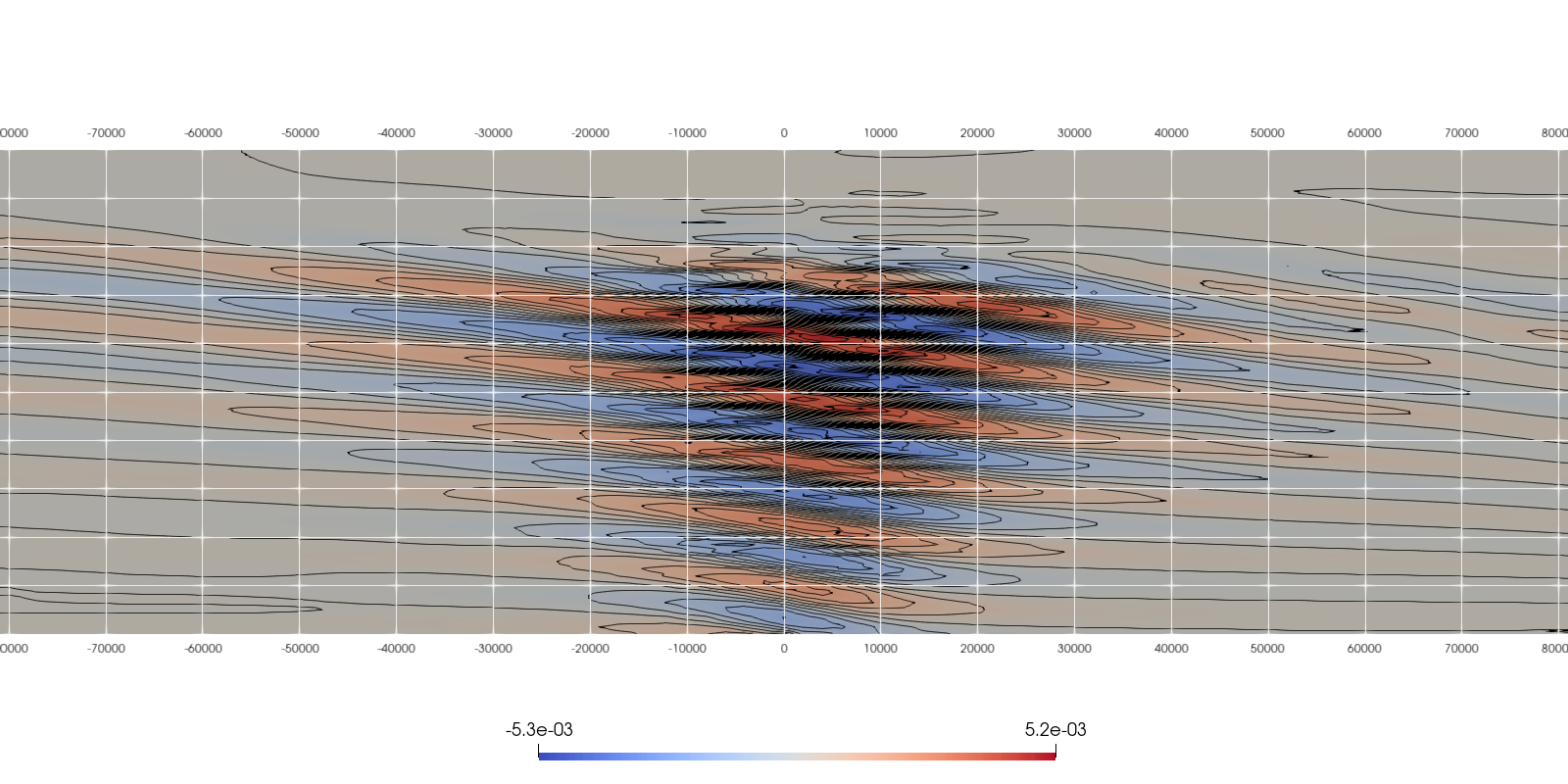}}
  \caption{\label{fig:ag hy}Contours of
  vertical velocity for the hydrostatic mountain wave test at
  $t=1.5\times 10^4s$. Contour intervals are every $5\times
  10^{-4}ms^{-1}$\response{; horizontal lines are plotted every 5000m.}}
\end{figure}

Contour plots of the vertical velocity are provided in Figures
\ref{fig:ag nh} and \ref{fig:ag hy}, respectively. We observe good
agreement with the solutions plotted in \citet{melvin2010inherently}.

\subsection{Sch\"ar test}
\citet{schar2002new} describe a more challenging mountain
wave test with a mountain range orography that varies over multiple length scales, defined as
\begin{equation}
  z_s(x) = h_m\exp\left(-\left(\frac{x}{a}\right)^2\right)
  \cos^2\left(\frac{\pi x}{\lambda}\right),
\end{equation}
where $h_m=250m$, $\lambda=4\times 10^3m$, and $a=5\times 10^3m$. The
domain is given by $-L/2\leq x \leq L/2$ where $L=10^5m$ and $0\leq z
\leq H = 3\times 10^4m$.  The initial stratification is initialised in
the same manner as for the nonhydrostatic mountain wave and we also
add an absorbing term as described in equation \ref{eq:absorb}, with
$z_B=H-10^4=2\times10^4m$ and $\bar{\mu}\Delta t=1.2s$.  The velocity
is initially horizontal with $u=10ms^{-1}$.  Following
\citet{melvin2010inherently}, we ran the test with two values of
$\Delta t$, $8s$ and $40s$. As might be expected from a solution that
has almost reached a steady state for the lee wave pattern and the
fact that we are not using a splitting method, the numerical solutions
obtained were indistinguishable. The larger timestep requires more
linear solver iterations, as displayed in Table \ref{tab:stats}.  This
is offset by the larger timestep, so the time to solution is similar
(shorter for the larger timestep, in fact). This reflects the fact
that our columnwise preconditioner produces mesh independent
convergence rates when applied to the linearisation about the state of
rest, but has a dependency on Courant number in general. The solutions
fit within the range of results obtained by others
e.g. \citet{straka1993numerical,giraldo2008study,bendall2020compatible},
given the turbulent nature of the solution.

\begin{figure}
  \centerline{\includegraphics[width=18cm]{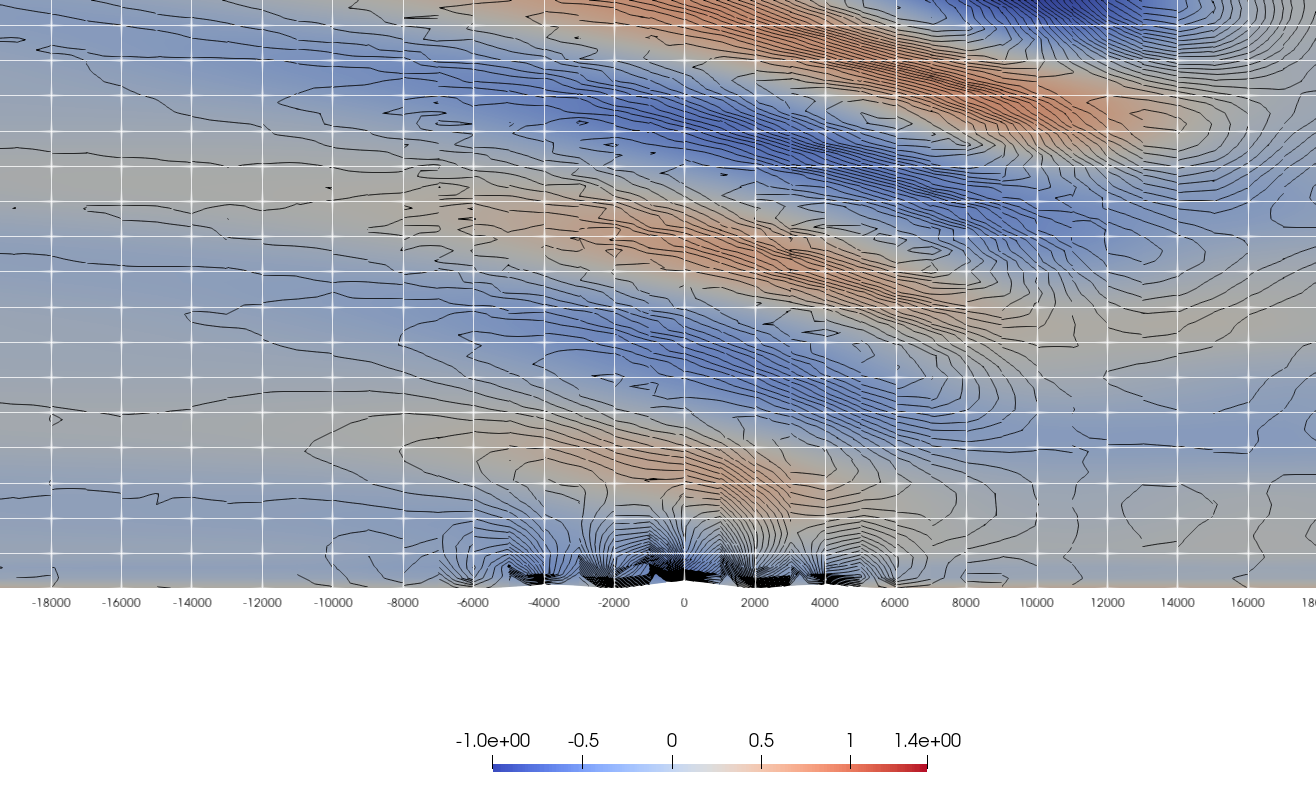}}
  \caption{\label{fig:schaer}Contours of vertical velocity for the
    Sch\"ar mountain wave test at $t=1.8\times 10^4s$, with $\Delta
    t=8s$. Contour intervals are every $5\times 10^{-2}ms^{-1}$.}
\end{figure}

\section{Summary and outlook}
\label{sec:summary}

In this paper we presented a compatible finite element discretisation
for the compressible Euler equations, and demonstrated numerical
robustness using a standard suite of vertical slice tests for
numerical weather prediction. In all cases the results are very
similar to published results. In future work we will demonstrate the
discretisation in the fully three dimensional setting on the sphere.

One novel feature of our approach is a monolithic approach to solving
a fully implicit system. This approach shows promise for producing
robust timestepping approaches. In experiments with the shallow water
equations we have used the additive Schwarz approach as a smoother for
a multigrid scheme, which has led to faster convergence of the
iterative solver. We were unable to do that in this work because
Firedrake does not currently support mesh hierarchies on periodic
meshes which are required for this suite of tests. We will rectify
this in further work.

One reason for our interest in fully implicit
methods is that these are required for the time parallel algorithms
that we are currently developing. It is also interesting to consider
Rosenbrock methods that only require the solution of linear systems
(as opposed to the nonlinear systems coming from the implicit midpoint
rule) linearised about the state at the start of the timestep. In
particular, we plan to experiment with a Rosenbrock version of the
TR-BDF2 timestepping scheme, previously applied to numerical weather
prediction in \citet{tumolo2015semi}.

\section*{Funding and Conflicts of Interest}
The authors would like to acknowledge funding from UKRI and The Met
Office (SPF ExCALIBUR EX20-8), the EPSRC PRISM platform grants
EP/L000407/1 and EP/R029423/1, EPSRC grant EP/R029628/1 and NERC
grants NE/I02013X/1, NE/R008795/1 and NE/M013634/1.

The authors have no relevant financial or non-financial interests to
disclose.

\bibliography{cfm_slice}

\begin{thebibliography}{48}
\expandafter\ifx\csname natexlab\endcsname\relax\def\natexlab#1{#1}\fi
\providecommand{\url}[1]{\texttt{#1}}
\providecommand{\href}[2]{#2}
\providecommand{\path}[1]{#1}
\providecommand{\DOIprefix}{doi:}
\providecommand{\ArXivprefix}{arXiv:}
\providecommand{\URLprefix}{URL: }
\providecommand{\Pubmedprefix}{pmid:}
\providecommand{\doi}[1]{\href{http://dx.doi.org/#1}{\path{#1}}}
\providecommand{\Pubmed}[1]{\href{pmid:#1}{\path{#1}}}
\providecommand{\bibinfo}[2]{#2}
\ifx\xfnm\relax \def\xfnm[#1]{\unskip,\space#1}\fi
\bibitem[{Adler et~al.(2021)Adler, Benson, Cyr, Farrell, MacLachlan and
  Tuminaro}]{adler2021monolithic}
\bibinfo{author}{Adler, J.H.}, \bibinfo{author}{Benson, T.R.},
  \bibinfo{author}{Cyr, E.C.}, \bibinfo{author}{Farrell, P.E.},
  \bibinfo{author}{MacLachlan, S.P.}, \bibinfo{author}{Tuminaro, R.S.},
  \bibinfo{year}{2021}.
\newblock \bibinfo{title}{Monolithic multigrid methods for
  magnetohydrodynamics}.
\newblock \bibinfo{journal}{SIAM Journal on Scientific Computing}
  \bibinfo{volume}{43}, \bibinfo{pages}{S70--S91}.
\bibitem[{Aln{\ae}s et~al.(2014)Aln{\ae}s, Logg, {\O}lgaard, Rognes and
  Wells}]{alnaes2014unified}
\bibinfo{author}{Aln{\ae}s, M.S.}, \bibinfo{author}{Logg, A.},
  \bibinfo{author}{{\O}lgaard, K.B.}, \bibinfo{author}{Rognes, M.E.},
  \bibinfo{author}{Wells, G.N.}, \bibinfo{year}{2014}.
\newblock \bibinfo{title}{Unified form language: A domain-specific language for
  weak formulations of partial differential equations}.
\newblock \bibinfo{journal}{ACM Transactions on Mathematical Software (TOMS)}
  \bibinfo{volume}{40}, \bibinfo{pages}{1--37}.
\bibitem[{Arakawa and Lamb(1981)}]{arakawa1981potential}
\bibinfo{author}{Arakawa, A.}, \bibinfo{author}{Lamb, V.R.},
  \bibinfo{year}{1981}.
\newblock \bibinfo{title}{A potential enstrophy and energy conserving scheme
  for the shallow water equations}.
\newblock \bibinfo{journal}{Monthly Weather Review} \bibinfo{volume}{109},
  \bibinfo{pages}{18--36}.
\bibitem[{Balay et~al.(2020)Balay, Abhyankar, Adams, Brown, Brune, Buschelman,
  Dalcin, Dener, Eijkhout, Gropp et~al.}]{balay2020petsc}
\bibinfo{author}{Balay, S.}, \bibinfo{author}{Abhyankar, S.},
  \bibinfo{author}{Adams, M.F.}, \bibinfo{author}{Brown, J.},
  \bibinfo{author}{Brune, P.}, \bibinfo{author}{Buschelman, K.},
  \bibinfo{author}{Dalcin, L.}, \bibinfo{author}{Dener, A.},
  \bibinfo{author}{Eijkhout, V.}, \bibinfo{author}{Gropp, W.}, et~al.,
  \bibinfo{year}{2020}.
\newblock \bibinfo{title}{PETSc Users Manual (Rev. 3.13)}.
\newblock \bibinfo{type}{Technical Report}. Argonne National Lab.(ANL),
  Argonne, IL (United States).
\bibitem[{Bauer and Cotter(2018)}]{bauer2018energy}
\bibinfo{author}{Bauer, W.}, \bibinfo{author}{Cotter, C.J.},
  \bibinfo{year}{2018}.
\newblock \bibinfo{title}{Energy--enstrophy conserving compatible finite
  element schemes for the rotating shallow water equations with slip boundary
  conditions}.
\newblock \bibinfo{journal}{Journal of Computational Physics}
  \bibinfo{volume}{373}, \bibinfo{pages}{171--187}.
\bibitem[{Bendall et~al.(2020)Bendall, Gibson, Shipton, Cotter and
  Shipway}]{bendall2020compatible}
\bibinfo{author}{Bendall, T.M.}, \bibinfo{author}{Gibson, T.H.},
  \bibinfo{author}{Shipton, J.}, \bibinfo{author}{Cotter, C.J.},
  \bibinfo{author}{Shipway, B.}, \bibinfo{year}{2020}.
\newblock \bibinfo{title}{A compatible finite-element discretisation for the
  moist compressible {E}uler equations}.
\newblock \bibinfo{journal}{Quarterly Journal of the Royal Meteorological
  Society} \bibinfo{volume}{146}, \bibinfo{pages}{3187--3205}.
\bibitem[{Bercea et~al.(2016)Bercea, McRae, Ham, Mitchell, Rathgeber, Nardi,
  Luporini and Kelly}]{bercea2016structure}
\bibinfo{author}{Bercea, G.T.}, \bibinfo{author}{McRae, A.T.},
  \bibinfo{author}{Ham, D.A.}, \bibinfo{author}{Mitchell, L.},
  \bibinfo{author}{Rathgeber, F.}, \bibinfo{author}{Nardi, L.},
  \bibinfo{author}{Luporini, F.}, \bibinfo{author}{Kelly, P.H.},
  \bibinfo{year}{2016}.
\newblock \bibinfo{title}{A structure-exploiting numbering algorithm for finite
  elements on extruded meshes, and its performance evaluation in {F}iredrake}.
\newblock \bibinfo{journal}{Geoscientific Model Development}
  \bibinfo{volume}{9}, \bibinfo{pages}{3803--3815}.
\bibitem[{Betteridge et~al.(2023)Betteridge, Cotter, Gibson, Griffith, Melvin
  and M{\"u}ller}]{betteridge2022hybridised}
\bibinfo{author}{Betteridge, J.D.}, \bibinfo{author}{Cotter, C.J.},
  \bibinfo{author}{Gibson, T.H.}, \bibinfo{author}{Griffith, M.J.},
  \bibinfo{author}{Melvin, T.}, \bibinfo{author}{M{\"u}ller, E.H.},
  \bibinfo{year}{2023}.
\newblock \bibinfo{title}{Hybridised multigrid preconditioners for a compatible
  finite element dynamical core}.
\newblock \bibinfo{journal}{to appear in Quarterly Journal of the Royal
  Meteorological Society} .
\bibitem[{Burman(2005)}]{burman2005unified}
\bibinfo{author}{Burman, E.}, \bibinfo{year}{2005}.
\newblock \bibinfo{title}{A unified analysis for conforming and nonconforming
  stabilized finite element methods using interior penalty}.
\newblock \bibinfo{journal}{SIAM journal on numerical analysis}
  \bibinfo{volume}{43}, \bibinfo{pages}{2012--2033}.
\bibitem[{Burman and Ern(2007)}]{burman2007continuous}
\bibinfo{author}{Burman, E.}, \bibinfo{author}{Ern, A.}, \bibinfo{year}{2007}.
\newblock \bibinfo{title}{Continuous interior penalty hp-finite element methods
  for advection and advection-diffusion equations}.
\newblock \bibinfo{journal}{Mathematics of computation} \bibinfo{volume}{76},
  \bibinfo{pages}{1119--1140}.
\bibitem[{Cockburn et~al.(2007)Cockburn, Kanschat and
  Sch{\"o}tzau}]{cockburn2007note}
\bibinfo{author}{Cockburn, B.}, \bibinfo{author}{Kanschat, G.},
  \bibinfo{author}{Sch{\"o}tzau, D.}, \bibinfo{year}{2007}.
\newblock \bibinfo{title}{A note on discontinuous {G}alerkin divergence-free
  solutions of the {N}avier--{S}tokes equations}.
\newblock \bibinfo{journal}{Journal of Scientific Computing}
  \bibinfo{volume}{31}, \bibinfo{pages}{61--73}.
\bibitem[{Cotter and Shipton(2012)}]{cotter2012mixed}
\bibinfo{author}{Cotter, C.J.}, \bibinfo{author}{Shipton, J.},
  \bibinfo{year}{2012}.
\newblock \bibinfo{title}{Mixed finite elements for numerical weather
  prediction}.
\newblock \bibinfo{journal}{Journal of Computational Physics}
  \bibinfo{volume}{231}, \bibinfo{pages}{7076--7091}.
\bibitem[{Cotter and Thuburn(2014)}]{cotter2014finite}
\bibinfo{author}{Cotter, C.J.}, \bibinfo{author}{Thuburn, J.},
  \bibinfo{year}{2014}.
\newblock \bibinfo{title}{A finite element exterior calculus framework for the
  rotating shallow-water equations}.
\newblock \bibinfo{journal}{Journal of Computational Physics}
  \bibinfo{volume}{257}, \bibinfo{pages}{1506--1526}.
\bibitem[{Danilov(2010)}]{danilov2010utility}
\bibinfo{author}{Danilov, S.}, \bibinfo{year}{2010}.
\newblock \bibinfo{title}{On utility of triangular {C}-grid type discretization
  for numerical modeling of large-scale ocean flows}.
\newblock \bibinfo{journal}{Ocean Dynamics} \bibinfo{volume}{60},
  \bibinfo{pages}{1361--1369}.
\bibitem[{Dubos et~al.(2015)Dubos, Dubey, Tort, Mittal, Meurdesoif and
  Hourdin}]{dubos2015dynamico}
\bibinfo{author}{Dubos, T.}, \bibinfo{author}{Dubey, S.},
  \bibinfo{author}{Tort, M.}, \bibinfo{author}{Mittal, R.},
  \bibinfo{author}{Meurdesoif, Y.}, \bibinfo{author}{Hourdin, F.},
  \bibinfo{year}{2015}.
\newblock \bibinfo{title}{{DYNAMICO}-1.0, an icosahedral hydrostatic dynamical
  core designed for consistency and versatility}.
\newblock \bibinfo{journal}{Geoscientific Model Development}
  \bibinfo{volume}{8}, \bibinfo{pages}{3131--3150}.
\bibitem[{Eldred et~al.(2019)Eldred, Dubos and Kritsikis}]{eldred2019quasi}
\bibinfo{author}{Eldred, C.}, \bibinfo{author}{Dubos, T.},
  \bibinfo{author}{Kritsikis, E.}, \bibinfo{year}{2019}.
\newblock \bibinfo{title}{A quasi-{H}amiltonian discretization of the thermal
  shallow water equations}.
\newblock \bibinfo{journal}{Journal of Computational Physics}
  \bibinfo{volume}{379}, \bibinfo{pages}{1--31}.
\bibitem[{Farrell et~al.(2019)Farrell, Mitchell and
  Wechsung}]{farrell2019augmented}
\bibinfo{author}{Farrell, P.E.}, \bibinfo{author}{Mitchell, L.},
  \bibinfo{author}{Wechsung, F.}, \bibinfo{year}{2019}.
\newblock \bibinfo{title}{An augmented {Lagrangian preconditioner for the 3D
  stationary incompressible Navier--Stokes equations at high Reynolds number}}.
\newblock \bibinfo{journal}{SIAM Journal on Scientific Computing}
  \bibinfo{volume}{41}, \bibinfo{pages}{A3073--A3096}.
\bibitem[{Gassmann(2013)}]{gassmann2013global}
\bibinfo{author}{Gassmann, A.}, \bibinfo{year}{2013}.
\newblock \bibinfo{title}{A global hexagonal {C}-grid non-hydrostatic dynamical
  core ({ICON-IAP}) designed for energetic consistency}.
\newblock \bibinfo{journal}{Quarterly Journal of the Royal Meteorological
  Society} \bibinfo{volume}{139}, \bibinfo{pages}{152--175}.
\bibitem[{Gawlik and Gay-Balmaz(2020)}]{gawlik2020conservative}
\bibinfo{author}{Gawlik, E.S.}, \bibinfo{author}{Gay-Balmaz, F.},
  \bibinfo{year}{2020}.
\newblock \bibinfo{title}{A conservative finite element method for the
  incompressible {E}uler equations with variable density}.
\newblock \bibinfo{journal}{Journal of Computational Physics}
  \bibinfo{volume}{412}, \bibinfo{pages}{109439}.
\bibitem[{Gawlik and Gay-Balmaz(2022)}]{gawlik2022finite}
\bibinfo{author}{Gawlik, E.S.}, \bibinfo{author}{Gay-Balmaz, F.},
  \bibinfo{year}{2022}.
\newblock \bibinfo{title}{A finite element method for {MHD} that preserves
  energy, cross-helicity, magnetic helicity, incompressibility, and {div B=
  0}}.
\newblock \bibinfo{journal}{Journal of Computational Physics}
  \bibinfo{volume}{450}, \bibinfo{pages}{110847}.
\bibitem[{Giraldo and Restelli(2008)}]{giraldo2008study}
\bibinfo{author}{Giraldo, F.X.}, \bibinfo{author}{Restelli, M.},
  \bibinfo{year}{2008}.
\newblock \bibinfo{title}{A study of spectral element and discontinuous
  {G}alerkin methods for the {N}avier--{S}tokes equations in nonhydrostatic
  mesoscale atmospheric modeling: Equation sets and test cases}.
\newblock \bibinfo{journal}{Journal of Computational Physics}
  \bibinfo{volume}{227}, \bibinfo{pages}{3849--3877}.
\bibitem[{Guerra and Ullrich(2016)}]{guerra2016high}
\bibinfo{author}{Guerra, J.E.}, \bibinfo{author}{Ullrich, P.A.},
  \bibinfo{year}{2016}.
\newblock \bibinfo{title}{A high-order staggered finite-element vertical
  discretization for non-hydrostatic atmospheric models}.
\newblock \bibinfo{journal}{Geoscientific Model Development}
  \bibinfo{volume}{9}, \bibinfo{pages}{2007--2029}.
\bibitem[{Laakmann et~al.(2022)Laakmann, Farrell and
  Mitchell}]{laakmann2022augmented}
\bibinfo{author}{Laakmann, F.}, \bibinfo{author}{Farrell, P.E.},
  \bibinfo{author}{Mitchell, L.}, \bibinfo{year}{2022}.
\newblock \bibinfo{title}{An augmented {L}agrangian preconditioner for the
  magnetohydrodynamics equations at high {R}eynolds and coupling numbers}.
\newblock \bibinfo{journal}{SIAM Journal on Scientific Computing}
  \bibinfo{volume}{44}, \bibinfo{pages}{B1018--B1044}.
\bibitem[{Lee and Palha(2018)}]{lee2018mixed}
\bibinfo{author}{Lee, D.}, \bibinfo{author}{Palha, A.}, \bibinfo{year}{2018}.
\newblock \bibinfo{title}{A mixed mimetic spectral element model of the
  rotating shallow water equations on the cubed sphere}.
\newblock \bibinfo{journal}{Journal of Computational Physics}
  \bibinfo{volume}{375}, \bibinfo{pages}{240--262}.
\bibitem[{Lee and Palha(2020)}]{lee2020mixed}
\bibinfo{author}{Lee, D.}, \bibinfo{author}{Palha, A.}, \bibinfo{year}{2020}.
\newblock \bibinfo{title}{A mixed mimetic spectral element model of the 3d
  compressible {E}uler equations on the cubed sphere}.
\newblock \bibinfo{journal}{Journal of Computational Physics}
  \bibinfo{volume}{401}, \bibinfo{pages}{108993}.
\bibitem[{McRae et~al.(2016)McRae, Bercea, Mitchell, Ham and
  Cotter}]{mcrae2016automated}
\bibinfo{author}{McRae, A.T.}, \bibinfo{author}{Bercea, G.T.},
  \bibinfo{author}{Mitchell, L.}, \bibinfo{author}{Ham, D.A.},
  \bibinfo{author}{Cotter, C.J.}, \bibinfo{year}{2016}.
\newblock \bibinfo{title}{Automated generation and symbolic manipulation of
  tensor product finite elements}.
\newblock \bibinfo{journal}{SIAM Journal on Scientific Computing}
  \bibinfo{volume}{38}, \bibinfo{pages}{S25--S47}.
\bibitem[{McRae and Cotter(2014)}]{mcrae2014energy}
\bibinfo{author}{McRae, A.T.}, \bibinfo{author}{Cotter, C.J.},
  \bibinfo{year}{2014}.
\newblock \bibinfo{title}{Energy-and enstrophy-conserving schemes for the
  shallow-water equations, based on mimetic finite elements}.
\newblock \bibinfo{journal}{Quarterly Journal of the Royal Meteorological
  Society} \bibinfo{volume}{140}, \bibinfo{pages}{2223--2234}.
\bibitem[{Melvin et~al.(2019)Melvin, Benacchio, Shipway, Wood, Thuburn and
  Cotter}]{melvin2019mixed}
\bibinfo{author}{Melvin, T.}, \bibinfo{author}{Benacchio, T.},
  \bibinfo{author}{Shipway, B.}, \bibinfo{author}{Wood, N.},
  \bibinfo{author}{Thuburn, J.}, \bibinfo{author}{Cotter, C.},
  \bibinfo{year}{2019}.
\newblock \bibinfo{title}{A mixed finite-element, finite-volume, semi-implicit
  discretization for atmospheric dynamics: {C}artesian geometry}.
\newblock \bibinfo{journal}{Quarterly Journal of the Royal Meteorological
  Society} \bibinfo{volume}{145}, \bibinfo{pages}{2835--2853}.
\bibitem[{Melvin et~al.(2010)Melvin, Dubal, Wood, Staniforth and
  Zerroukat}]{melvin2010inherently}
\bibinfo{author}{Melvin, T.}, \bibinfo{author}{Dubal, M.},
  \bibinfo{author}{Wood, N.}, \bibinfo{author}{Staniforth, A.},
  \bibinfo{author}{Zerroukat, M.}, \bibinfo{year}{2010}.
\newblock \bibinfo{title}{An inherently mass-conserving iterative semi-implicit
  semi-{L}agrangian discretization of the non-hydrostatic vertical-slice
  equations}.
\newblock \bibinfo{journal}{Quarterly Journal of the Royal Meteorological
  Society: A journal of the atmospheric sciences, applied meteorology and
  physical oceanography} \bibinfo{volume}{136}, \bibinfo{pages}{799--814}.
\bibitem[{Natale and Cotter(2017)}]{natale2017scale}
\bibinfo{author}{Natale, A.}, \bibinfo{author}{Cotter, C.J.},
  \bibinfo{year}{2017}.
\newblock \bibinfo{title}{Scale-selective dissipation in energy-conserving
  finite-element schemes for two-dimensional turbulence}.
\newblock \bibinfo{journal}{Quarterly Journal of the Royal Meteorological
  Society} \bibinfo{volume}{143}, \bibinfo{pages}{1734--1745}.
\bibitem[{Natale and Cotter(2018)}]{natale2018variational}
\bibinfo{author}{Natale, A.}, \bibinfo{author}{Cotter, C.J.},
  \bibinfo{year}{2018}.
\newblock \bibinfo{title}{A variational finite-element discretization approach
  for perfect incompressible fluids}.
\newblock \bibinfo{journal}{IMA Journal of Numerical Analysis}
  \bibinfo{volume}{38}, \bibinfo{pages}{1388--1419}.
\bibitem[{Natale et~al.(2016)Natale, Shipton and Cotter}]{natale2016compatible}
\bibinfo{author}{Natale, A.}, \bibinfo{author}{Shipton, J.},
  \bibinfo{author}{Cotter, C.J.}, \bibinfo{year}{2016}.
\newblock \bibinfo{title}{Compatible finite element spaces for geophysical
  fluid dynamics}.
\newblock \bibinfo{journal}{Dynamics and Statistics of the Climate System}
  \bibinfo{volume}{1}.
\bibitem[{Rathgeber et~al.(2016)Rathgeber, Ham, Mitchell, Lange, Luporini,
  McRae, Bercea, Markall and Kelly}]{rathgeber2016firedrake}
\bibinfo{author}{Rathgeber, F.}, \bibinfo{author}{Ham, D.A.},
  \bibinfo{author}{Mitchell, L.}, \bibinfo{author}{Lange, M.},
  \bibinfo{author}{Luporini, F.}, \bibinfo{author}{McRae, A.T.},
  \bibinfo{author}{Bercea, G.T.}, \bibinfo{author}{Markall, G.R.},
  \bibinfo{author}{Kelly, P.H.}, \bibinfo{year}{2016}.
\newblock \bibinfo{title}{{F}iredrake: automating the finite element method by
  composing abstractions}.
\newblock \bibinfo{journal}{ACM Transactions on Mathematical Software (TOMS)}
  \bibinfo{volume}{43}, \bibinfo{pages}{1--27}.
\bibitem[{Rognes et~al.(2010)Rognes, Kirby and Logg}]{rognes2010efficient}
\bibinfo{author}{Rognes, M.E.}, \bibinfo{author}{Kirby, R.C.},
  \bibinfo{author}{Logg, A.}, \bibinfo{year}{2010}.
\newblock \bibinfo{title}{Efficient assembly of {H}(div) and {H}(curl)
  conforming finite elements}.
\newblock \bibinfo{journal}{SIAM Journal on Scientific Computing}
  \bibinfo{volume}{31}, \bibinfo{pages}{4130--4151}.
\bibitem[{Sch{\"a}r et~al.(2002)Sch{\"a}r, Leuenberger, Fuhrer, L{\"u}thi and
  Girard}]{schar2002new}
\bibinfo{author}{Sch{\"a}r, C.}, \bibinfo{author}{Leuenberger, D.},
  \bibinfo{author}{Fuhrer, O.}, \bibinfo{author}{L{\"u}thi, D.},
  \bibinfo{author}{Girard, C.}, \bibinfo{year}{2002}.
\newblock \bibinfo{title}{A new terrain-following vertical coordinate
  formulation for atmospheric prediction models}.
\newblock \bibinfo{journal}{Monthly Weather Review} \bibinfo{volume}{130},
  \bibinfo{pages}{2459--2480}.
\bibitem[{Sergeev et~al.(2023)Sergeev, Mayne, Bendall, Boutle, Brown, Kavcic,
  Kent, Kohary, Manners, Melvin et~al.}]{sergeev2023simulations}
\bibinfo{author}{Sergeev, D.E.}, \bibinfo{author}{Mayne, N.J.},
  \bibinfo{author}{Bendall, T.}, \bibinfo{author}{Boutle, I.A.},
  \bibinfo{author}{Brown, A.}, \bibinfo{author}{Kavcic, I.},
  \bibinfo{author}{Kent, J.}, \bibinfo{author}{Kohary, K.},
  \bibinfo{author}{Manners, J.}, \bibinfo{author}{Melvin, T.}, et~al.,
  \bibinfo{year}{2023}.
\newblock \bibinfo{title}{Simulations of idealised {3D} atmospheric flows on
  terrestrial planets using {LFRic-Atmosphere}}.
\newblock \bibinfo{journal}{arXiv preprint arXiv:2306.03614} .
\bibitem[{Shipton et~al.(2018)Shipton, Gibson and Cotter}]{shipton2018higher}
\bibinfo{author}{Shipton, J.}, \bibinfo{author}{Gibson, T.H.},
  \bibinfo{author}{Cotter, C.J.}, \bibinfo{year}{2018}.
\newblock \bibinfo{title}{Higher-order compatible finite element schemes for
  the nonlinear rotating shallow water equations on the sphere}.
\newblock \bibinfo{journal}{Journal of Computational Physics}
  \bibinfo{volume}{375}, \bibinfo{pages}{1121--1137}.
\bibitem[{Skamarock and Klemp(1994)}]{skamarock1994efficiency}
\bibinfo{author}{Skamarock, W.C.}, \bibinfo{author}{Klemp, J.B.},
  \bibinfo{year}{1994}.
\newblock \bibinfo{title}{Efficiency and accuracy of the {K}lemp-{W}ilhelmson
  time-splitting technique}.
\newblock \bibinfo{journal}{Monthly Weather Review} \bibinfo{volume}{122},
  \bibinfo{pages}{2623--2630}.
\bibitem[{Staniforth and Thuburn(2012)}]{staniforth2012horizontal}
\bibinfo{author}{Staniforth, A.}, \bibinfo{author}{Thuburn, J.},
  \bibinfo{year}{2012}.
\newblock \bibinfo{title}{Horizontal grids for global weather and climate
  prediction models: a review}.
\newblock \bibinfo{journal}{Quarterly Journal of the Royal Meteorological
  Society} \bibinfo{volume}{138}, \bibinfo{pages}{1--26}.
\bibitem[{Straka et~al.(1993)Straka, Wilhelmson, Wicker, Anderson and
  Droegemeier}]{straka1993numerical}
\bibinfo{author}{Straka, J.M.}, \bibinfo{author}{Wilhelmson, R.B.},
  \bibinfo{author}{Wicker, L.J.}, \bibinfo{author}{Anderson, J.R.},
  \bibinfo{author}{Droegemeier, K.K.}, \bibinfo{year}{1993}.
\newblock \bibinfo{title}{Numerical solutions of a non-linear density current:
  A benchmark solution and comparisons}.
\newblock \bibinfo{journal}{International Journal for Numerical Methods in
  Fluids} \bibinfo{volume}{17}, \bibinfo{pages}{1--22}.
\bibitem[{Taylor et~al.(2020)Taylor, Guba, Steyer, Ullrich, Hall and
  Eldred}]{taylor2020energy}
\bibinfo{author}{Taylor, M.A.}, \bibinfo{author}{Guba, O.},
  \bibinfo{author}{Steyer, A.}, \bibinfo{author}{Ullrich, P.A.},
  \bibinfo{author}{Hall, D.M.}, \bibinfo{author}{Eldred, C.},
  \bibinfo{year}{2020}.
\newblock \bibinfo{title}{An energy consistent discretization of the
  nonhydrostatic equations in primitive variables}.
\newblock \bibinfo{journal}{Journal of Advances in Modeling Earth Systems}
  \bibinfo{volume}{12}, \bibinfo{pages}{e2019MS001783}.
\bibitem[{Thuburn et~al.(2014)Thuburn, Cotter and Dubos}]{thuburn2014mimetic}
\bibinfo{author}{Thuburn, J.}, \bibinfo{author}{Cotter, C.},
  \bibinfo{author}{Dubos, T.}, \bibinfo{year}{2014}.
\newblock \bibinfo{title}{A mimetic, semi-implicit, forward-in-time, finite
  volume shallow water model: comparison of hexagonal--icosahedral and
  cubed-sphere grids}.
\newblock \bibinfo{journal}{Geoscientific Model Development}
  \bibinfo{volume}{7}, \bibinfo{pages}{909--929}.
\bibitem[{Thuburn and Cotter(2012)}]{thuburn2012framework}
\bibinfo{author}{Thuburn, J.}, \bibinfo{author}{Cotter, C.J.},
  \bibinfo{year}{2012}.
\newblock \bibinfo{title}{A framework for mimetic discretization of the
  rotating shallow-water equations on arbitrary polygonal grids}.
\newblock \bibinfo{journal}{SIAM Journal on Scientific Computing}
  \bibinfo{volume}{34}, \bibinfo{pages}{B203--B225}.
\bibitem[{Tumolo and Bonaventura(2015)}]{tumolo2015semi}
\bibinfo{author}{Tumolo, G.}, \bibinfo{author}{Bonaventura, L.},
  \bibinfo{year}{2015}.
\newblock \bibinfo{title}{A semi-implicit, semi-{L}agrangian discontinuous
  {G}alerkin framework for adaptive numerical weather prediction}.
\newblock \bibinfo{journal}{Quarterly Journal of the Royal Meteorological
  Society} \bibinfo{volume}{141}, \bibinfo{pages}{2582--2601}.
\bibitem[{Wimmer et~al.(2020)Wimmer, Cotter and Bauer}]{wimmer2020energy}
\bibinfo{author}{Wimmer, G.A.}, \bibinfo{author}{Cotter, C.J.},
  \bibinfo{author}{Bauer, W.}, \bibinfo{year}{2020}.
\newblock \bibinfo{title}{Energy conserving upwinded compatible finite element
  schemes for the rotating shallow water equations}.
\newblock \bibinfo{journal}{Journal of Computational Physics}
  \bibinfo{volume}{401}, \bibinfo{pages}{109016}.
\bibitem[{Wimmer et~al.(2021)Wimmer, Cotter and Bauer}]{wimmer2021energy}
\bibinfo{author}{Wimmer, G.A.}, \bibinfo{author}{Cotter, C.J.},
  \bibinfo{author}{Bauer, W.}, \bibinfo{year}{2021}.
\newblock \bibinfo{title}{Energy conserving {SUPG} methods for compatible
  finite element schemes in numerical weather prediction}.
\newblock \bibinfo{journal}{The SMAI journal of computational mathematics}
  \bibinfo{volume}{7}, \bibinfo{pages}{267--300}.
\bibitem[{Wood et~al.(2014)Wood, Staniforth, White, Allen, Diamantakis, Gross,
  Melvin, Smith, Vosper, Zerroukat et~al.}]{wood2014inherently}
\bibinfo{author}{Wood, N.}, \bibinfo{author}{Staniforth, A.},
  \bibinfo{author}{White, A.}, \bibinfo{author}{Allen, T.},
  \bibinfo{author}{Diamantakis, M.}, \bibinfo{author}{Gross, M.},
  \bibinfo{author}{Melvin, T.}, \bibinfo{author}{Smith, C.},
  \bibinfo{author}{Vosper, S.}, \bibinfo{author}{Zerroukat, M.}, et~al.,
  \bibinfo{year}{2014}.
\newblock \bibinfo{title}{An inherently mass-conserving semi-implicit
  semi-{L}agrangian discretization of the deep-atmosphere global
  non-hydrostatic equations}.
\newblock \bibinfo{journal}{Quarterly Journal of the Royal Meteorological
  Society} \bibinfo{volume}{140}, \bibinfo{pages}{1505--1520}.
\bibitem[{Yamazaki et~al.(2017)Yamazaki, Shipton, Cullen, Mitchell and
  Cotter}]{yamazaki2017vertical}
\bibinfo{author}{Yamazaki, H.}, \bibinfo{author}{Shipton, J.},
  \bibinfo{author}{Cullen, M.J.}, \bibinfo{author}{Mitchell, L.},
  \bibinfo{author}{Cotter, C.J.}, \bibinfo{year}{2017}.
\newblock \bibinfo{title}{Vertical slice modelling of nonlinear {E}ady waves
  using a compatible finite element method}.
\newblock \bibinfo{journal}{Journal of Computational Physics}
  \bibinfo{volume}{343}, \bibinfo{pages}{130--149}.

\end{thebibliography}
\end{document}